\documentclass[parskip=half,bibliography=totoc]{scrartcl}

\usepackage[automark]{scrlayer-scrpage}								
\usepackage{amsfonts}												
\usepackage{amsmath}												
\usepackage{mathtools}												
\usepackage{stmaryrd}												
\usepackage{amssymb}												
\usepackage{extarrows}												
\usepackage{dsfont} 												
\usepackage{mathrsfs}												
\usepackage{accents}												
\usepackage[T1]{fontenc}											
\usepackage[utf8]{inputenc}											
\usepackage{subcaption} 											
\usepackage[top=1in,bottom=1in,left=1.25in,right=1.25in]{geometry}	
\usepackage{footmisc}												
\usepackage{enumerate} 												
\usepackage{booktabs}												
\usepackage{authblk}												
\usepackage{abstract}												
\usepackage{etoolbox}												

\usepackage{tikz-cd}												

\usepackage{amsthm}													

\usepackage[colorlinks=true,citecolor=blue]{hyperref}				
\usepackage[noabbrev]{cleveref}										
\crefname{lem}{lemma}{lemmata}

\newcommand{\pd}[2]{\frac{\partial #1}{\partial #2}}
\renewcommand{\d}{\mathrm{d}}

\newcommand{\CP}{\mathbb{C}\mathrm{P}}
\newcommand{\CH}{\mathbb{C} H}

\newcommand{\Z}{\mathbb{Z}}

\newcommand{\R}{\mathbb{R}}
\newcommand{\C}{\mathbb{C}}
\renewcommand{\H}{\mathbb{H}}

\renewcommand{\Im}{\operatorname{Im}}

\DeclareMathOperator{\Span}{span}

\DeclareMathOperator{\End}{End}
\DeclareMathOperator{\Aut}{Aut}
\DeclareMathOperator{\Sym}{Sym}
\newcommand{\aut}{\mathfrak{aut}}

\DeclareMathOperator{\id}{id}

\DeclareMathOperator{\Heis}{Heis}

\newcommand{\h}{\mathrm{H}}
\newcommand{\q}{\mathrm{Q}}

\newcommand{\abs}[1]{\lvert #1 \rvert}
\newcommand{\norm}[1]{\lVert #1 \rVert}

\newcommand{\mf}[1]{\mathfrak{#1}}
\newcommand{\mc}[1]{\mathcal{#1}}

\newtheoremstyle{mythm}
{}
{}
{\slshape}
{}
{\bfseries\sffamily}
{.}
{ }
{}
\newtheoremstyle{mydef}
{}
{}
{}
{}
{\bfseries\sffamily}
{.}
{ }
{}

\theoremstyle{mythm}
\newtheorem{thm}{Theorem}[section]
\newtheorem{prop}[thm]{Proposition}
\newtheorem{cor}[thm]{Corollary}
\newtheorem{lem}[thm]{Lemma}
\theoremstyle{mydef}
\newtheorem{mydef}[thm]{Definition}
\newtheorem{rem}[thm]{Remark}

\newenvironment{myproof}[1][\proofname]{
	\proof[\sffamily\upshape#1]
}{\endproof}

\newcommand{\proofclear}{\hfill \qedsymbol}

\clearscrheadfoot
\ihead[]{\headmark}
\ohead[]{\pagemark}
\deffootnote[1em]{0em}{1em}{%
	\textsuperscript{\thefootnotemark}%
}
\setfootnoterule{3em}

\newcommand\numberthis{\stepcounter{equation}\tag{\theequation}}

\newenvironment{numberedlist}{\begin{enumerate}[\upshape(i)]}{\end{enumerate}}

\apptocmd{\sloppy}{\hbadness 10000\relax}{}{}

\title{Symmetries of quaternionic Kähler manifolds with $S^1$-symmetry}
\author{V.\ Cortés}
\author{A.\ Saha}
\author{D.\ Thung}
\affil{\normalsize  Department of Mathematics \\
University of Hamburg\\
Bundesstra\ss e 55, D-20146 Hamburg, Germany}

\date{}

\begin{document}
\maketitle

\begin{abstract}
	We study symmetry properties of quaternionic K\"ahler manifolds obtained by the HK/QK correspondence. To any Lie algebra $\mathfrak{g}$ of infinitesimal automorphisms of the initial hyper-K\"ahler data we associate a central extension of $\mathfrak{g}$, acting by infinitesimal automorphisms of the resulting quaternionic K\"ahler manifold. More specifically, we study the metrics obtained by the one-loop deformation of the $c$-map construction, proving that the Lie algebra of infinitesimal automorphisms of the initial projective special K\"ahler manifold gives rise to a Lie algebra of Killing fields of the corresponding one-loop deformed $c$-map space. As an application, we show that this construction increases the cohomogeneity of the automorphism groups by at most one. In particular, if the initial manifold is homogeneous then the one-loop deformed metric is of cohomogeneity at most one. As an example, we consider the one-loop deformation of the symmetric quaternionic K\"ahler metric on $SU(n,2)/S(U(n)\times U(2))$, which we prove is of cohomogeneity exactly one. This family generalizes the so-called universal hypermultiplet ($n=1$), for which we determine the full isometry group.\par
\end{abstract}

\clearpage

\tableofcontents
\clearpage

\section{Introduction}

In the late 1980's, physicists working on string theory discovered a new method of constructing quaternionic K\"ahler manifolds of negative scalar curvature, which was dubbed the (supergravity) c-map. A local description of the resulting metrics was first given by Ferrara and Sabharwal~\cite{FS1990}, and they were described in differential-geometric terms by Hitchin~\cite{Hit2009}.

The input for the $c$-map (recalled in more detail in \Cref{sec:prelims}) is a projective special K\"ahler (PSK) manifold, which we will denote by $\bar M$. It arises as the K\"ahler quotient by a distinguished circle action on an associated manifold $M$, which caries a conical affine special K\"ahler (CASK) structure. As shown in~\cite{ACM2013,ACDM2015}, the supergravity $c$-map can be described as a two-step process. The first, called the rigid $c$-map, consists of passing from the corresponding CASK manifold $M$ to its cotangent bundle $N\coloneqq T^*M$, and defining a natural (pseudo-)hyper-K\"ahler structure on the latter. The circle action on $M$ lifts to a so-called rotating (circle) symmetry of the hyper-K\"ahler structure on $N$, and plays a crucial role in the the second step, the HK/QK correspondence. This correspondence produces a quaternionic K\"ahler manifold $\bar N$ of the same dimension out of any hyper-K\"ahler manifold with rotating symmetry, completing the following diagram:
\begin{equation}\label{eq:diagram}
	\begin{tikzcd}[column sep=2.5cm, row sep=1cm]
		M \ar[r,"\text{rigid }c"] \ar[d,"\C^*"'] & N \ar[d,"\text{HK/QK}"]\\
		\bar M \ar[r,"\text{supergravity }c"'] & \bar N
	\end{tikzcd}
\end{equation}
A one-parameter deformation of the $c$-map was discovered by Robles-Llana, Sauer\-essig and Vandoren \cite{RSV2006}. It is known as the one-loop deformed $c$-map because of its physical origin as a perturbative quantum correction to the Ferrara-Sabharwal metric. It is obtained by modifying the Ferrara--Sabharwal metric by introducing a real parameter $c$. When $c=0$, one recovers the Ferrara--Sabharwal metric. Curiously, $c$ has a very simple meaning on the hyper-K\"ahler side: It parametrizes the additive freedom in choosing a Hamiltonian for the lifted circle action. This key observation means that, by working on the hyper-K\"ahler side of the correspondence, we can treat the $c$-map and its one-loop correction simultaneously.

The HK/QK correspondence has been a subject of interest for both physicists and mathematicians (see, for instance,~\cite{APP2011,ACM2013,Hit2013,Hit2014,ACDM2015,MS2015}). It was first discovered by Haydys in the Riemannian setting~\cite{Hay2008} and subsequently generalized to allow for indefinite metrics~\cite{ACM2013,ACDM2015}. It is this latter form which can be used for the $c$-map construction. Nevertheless, the output of the $c$-map and its one-loop deformation is a positive-definite quaternionic K\"ahler manifold, which is moreover complete if the initial PSK manifold was complete and $c\geq 0$~\cite{CDS2017}. Summarizing and extending the previous remarks, the HK/QK correspondence assigns a quaternionic K\"ahler manifold with a nowhere-vanishing Killing field to any (possibly indefinite) hyper-K\"ahler manifold with rotating circle symmetry. Conversely, any quaternionic K\"ahler manifold with a Killing field locally arises in this fashion.

Another perspective on the HK/QK correspondence was advanced by Macia and Swann \cite{MS2015}: They showed that it can be recovered as a variation on the twist construction, a general duality between manifolds endowed with torus action, introduced by Swann in earlier work \cite{Swa2010}. One advantage of the twist construction, which we review in \Cref{sec:prelims} together with some facts about the relevant aspects of special K\"ahler geometry, is that it provides a precise prescription for carrying tensor fields on $N$ which are invariant under the natural circle action over to $\bar N$. Invariant vector fields, in particular, arise naturally in the context of the $c$-map. Infinitesimal automorphisms of the CASK structure on $M$ respect its distinguished circle action, and we show that they can be lifted to vector fields on $N$ which respect the rigid $c$-map structure on $N$.

With this in mind, we investigate the HK/QK correspondence in the presence of additional symmetries preserving the rigid $c$-map structure in \Cref{sec:symmetries}. We prove that, under the HK/QK correspondence as well as under its one-parameter deformation, these give rise to isometries of the resulting quaternionic K\"ahler manifold. Applying this to the $c$-map, we show that the (identity component of the) automorphism group of the PSK manifold $\bar M$ injects into the isometry group of its (one-loop deformed) $c$-map image $\bar N$. As a corollary, we find that if $\Aut \bar M$ acts with cohomogeneity $k$, then the isometry group of $\bar N$ acts with cohomogeneity at most $k+1$ for every $c\neq 0$, and at most $k$ for $c=0$ (the latter was previously known). 

In \Cref{sec:ex}, we apply our results to the study of a series of examples, namely the quaternionic K\"ahler manifolds obtained by applying the $c$-map to complex hyperbolic space $\CH^n=SU(n,1)/U(n)$, which is a homogeneous PSK manifold. The undeformed $c$-map yields the non-compact Wolf spaces $SU(n+1,2)/S(U(n+1)\times U(2))$. These can be characterized as the only complete, simply connected quaternionic K\"ahler manifolds with negative scalar curvature which are also K\"ahler. In the physics literature, the case $n=0$ is known as the ``universal hypermultiplet''. The one-loop deformed $c$-map yields a one-parameter deformation of the symmetric metric, and we investigate its isometry group. It follows from our results that it acts with cohomogeneity at most one and by using results from \cite{CST2020a}, where curvature formulas for the HK/QK correspondence are derived, we prove that the isometry group acts by cohomogeneity exactly one in all dimensions. This shows that our results are sharp in this family of examples and proves, in particular, that the one-loop deformed Ferrara--Sabharwal metric is not homogeneous, even locally. In the case of the universal hypermultiplet, we determine the full isometry group.

{\bfseries Acknowledgements}

D.T.~thanks Markus R\"oser for useful discussions. This work was supported by the German Science Foundation (DFG) under the Research Training Group 1670, and under Germany’s Excellence Strategy – EXC 2121 ``Quantum Universe'' – 390833306.

\clearpage

\section{Preliminaries}
\label{sec:prelims}

In this section, we set the stage for the rest of this paper, briefly recalling the most important background material. In particular, the terminology and notation introduced here will be used throughout \Cref{sec:symmetries}, where we prove our main results.

\subsection{Twist construction and the HK/QK correspondence}
\label{subsec:twist}

As mentioned in the introduction, Macia and Swann \cite{MS2015} studied the HK/QK correspondence from the point of view of Swann's twist construction. This is a general procedure which, given a manifold with torus action, produces a new manifold with twisted torus action, uniquely determined by a choice of so-called twist data on the original manifold. Moreover, geometric structures which are invariant under the torus action on the original manifold unambiguously correspond to invariant (twisted) structures on the twist manifold. The twist construction is a duality in the sense that the original manifold can be recovered by appropriately twisting the twist manifold.

In particular, the twist construction can be applied to a (pseudo-)hyper-K\"ahler manifold with circle action---as produced by the rigid $c$-map. Macia and Swann showed that, after deforming the original metric and choosing appropriate twist data, the twist of the deformed metric defines a quaternionic K\"ahler metric on the twist manifold. Their formulas reproduce the explicit expressions for the HK/QK correspondence found in \cite{ACDM2015}, casting it as a special case of the twist construction. In this section, we introduce the basic language of the twist construction and outline its application to the HK/QK correspondence, following \cite{MS2015}.

\subsubsection{The twist construction}

For our purposes, it suffices to describe the twist construction for circle actions; the interested reader is referred to \cite{Swa2010} for twists of higher-dimensional torus actions. 

Consider a manifold $N$ equipped with an $S^1$-action generated by a vector field $Z\in \mf X(N)$. To apply the twist construction, we consider an auxiliary manifold, namely the total space $P$ of an $S^1$-principal bundle $\pi_N:P\to N$. We furthermore equip this principal bundle with a principal $S^1$-connection $\eta\in \Omega^1(P)$ with curvature $\omega$. We wish to lift the  vector field $Z$ to $P$ so that the lift preserves the connection $\eta$ and commutes with the principal circle action. We will then obtain the twist of $M$ by dividing out the action generated by this lifted vector field.

\begin{prop}[\cite{Swa2010}]
	A lift as above exists if and only if $\iota_Z \omega$ is an exact one-form, i.e.\ $\iota_Z\omega=-\d f_Z$ for some Hamiltonian function $f_Z\in C^\infty(N)$.\proofclear 
\end{prop}

Concretely, the lift is given by $Z_P=\tilde Z+\pi_N^*f_Z X_P$, where $\tilde Z$ denotes the horizontal lift with respect to $\eta$, and $X_P$ generates the principal circle action. Note that one needs to fix a choice of Hamiltonian function. Thus, twisting involves the following twist data:
\begin{numberedlist}
	\item A vector field $Z\in \mf X(N)$ generating a circle action on $N$.
	\item A closed, integral two-form $\omega$ (i.e.~with integral periods) with respect to which $Z$ is Hamiltonian. Up to isomorphism, this form determines a unique principal circle bundle $\pi_N:P\to N$, which admits a connection $\eta$ whose curvature is $\omega$.
	\item A Hamiltonian function $f_Z$ for $Z$ with respect to $\omega$, which we assume is nowhere-vanishing (equivalently, the lifted vector field $Z_P$ is transverse to $\mc H\coloneqq \ker\eta$).
\end{numberedlist}

For given twist data, the twist of $N$ is defined as follows:

\begin{mydef}
	For $(Z,\omega,f_Z)$ as above, we call the quotient space $\bar N\coloneqq P/\langle Z_P\rangle$ the \emph{twist} of $N$ with respect to the twist data $(Z,\omega,f_Z)$.
\end{mydef}

If both $Z$ and $Z_P$ generate free and proper actions, the twist is a smooth manifold.

We now have the double fibration 
\begin{equation*}
	\begin{tikzcd}
		N & P \ar[l,"\pi_N"'] \ar[r,"\pi_{\bar N}"] & \bar N
	\end{tikzcd}
\end{equation*}
The horizontal distribution $\mc H\coloneqq \ker \eta$ is pointwise identified with the tangent spaces to $N$ as well as to $\bar N$, by means of the respective projections. This induces identifications between tensor fields in the following manner:

\begin{mydef}
	We say that vector fields $X\in \mf X(N)$ and $X'\in \mf X(\bar N)$ are \emph{$\mc H$-related}, and write $X'\sim_{\mc H} X$, if their horizontal lifts to $P$ (with respect to $\eta$) are equal. We say that differential forms $\alpha$ on $N$ and $\alpha'$ on $\bar N$ are $\mc H$-related if $(\pi_N^*\alpha)\circ P_{\mc H}=(\pi^*_{\bar N}\alpha')\circ P_{\mc H}$, where $P_{\mc H}:TP\to \mc H$, $U\mapsto U-\eta(U)X_P$, is the projection onto $\mc H$. The definition of $\mc H$-relatedness is now extended to arbitrary $(p,q)$-tensors by demanding compatibility with tensor products.
\end{mydef}

\begin{rem}\leavevmode
	\begin{numberedlist}
		\item Twisting is automatically compatible with contractions of tensor fields.
		\item The existence of a (well-defined) $\mc H$-related tensor field on $\bar N$ implies that the original tensor field on $N$ is $Z$-invariant. Thus, the twist construction only allows one to carry $Z$-invariant tensors on $N$ over to $\bar N$.
		\item The vector field $Z_{\bar N}\coloneqq (\pi_{\bar N})_*(X_P)\in \mf X(\bar N)$ plays a role dual to $Z$ in the following sense: If one performs a twist on $\bar N$ with respect to the data $(Z_{\bar N},(f_Z^{-1}\omega)',(f_Z^{-1})')$, where primes denote $\mc H$-relatedness, one recovers the original manifold $N$.
	\end{numberedlist}
\end{rem}

$Z$-invariant tensor fields on $N$ and, dually, $Z_{\bar N}$-invariant tensor fields on $\bar N$, are determined by their respective lifts to $P$. Therefore, one can uniquely determine $\mc H$-related tensor fields in terms of the twist data. For instance, the twist construction assigns to $X\in \mf X(N)$ the vector field $X'=(\pi_{\bar N})_*(\tilde X)$, where $\tilde X$ denotes the horizontal lift of $X$. One can also check how natural operations, such as the Lie bracket, behave under twisting:

\begin{lem}[\cite{Swa2010}]\leavevmode\label{lem:twisted-structures}
	\begin{numberedlist}
		\item If $h$ is a $Z$-invariant symmetric $(0,2)$-tensor on $N$, then the $\mc H$-related symmetric $(0,2)$-tensor $h'$ on $\bar N$ is determined by
		\begin{equation*}
			\pi_{\bar N}^*h'=\pi^*_Nh-2\eta\vee \pi_N^*(f_Z^{-1}\iota_Z h)
			+\eta^2\pi^*_N(f_Z^{-2}h(Z,Z))
		\end{equation*}
		\item If $\mf X(N)\ni X,Y\sim_{\mc H}X',Y'$, i.e.\ $\tilde X=\widehat{X'}$ and $\tilde Y=\widehat{Y'}$ (where the tilde and hat denote the respective horizontal lifts), we have 
		\begin{equation*}
			[X',Y']\sim_{\mc H}[X,Y]+f_Z^{-1}\omega(X,Y) Z
			\vspace{-1.1cm}
		\end{equation*}
	\end{numberedlist}
	\proofclear
\end{lem}

Since we will be interested in producing Killing fields on twist manifolds, we show explicitly how the Lie derivative of a symmetric $(0,2)$-tensor behaves under the twist:

\begin{lem}\label{lem:twisted-Killing}
	Let $V\sim_{\mc H} V'$ be vector fields and $h\sim_{\mc H} h'$ be symmetric $(0,2)$-tensors. Then $L_{V'}h'\sim_\mc{H} L_V h-2 f_Z^{-1} (\iota_V \omega)\vee (\iota_Z h)$.
\end{lem}
\begin{myproof}
	The Leibniz rule tells us that for two vector fields $X,Y$ on $M$, we have
	\begin{equation*}
		\pi_N^*\big((L_V h)(X,Y)\big)
		=L_{\tilde V}(\pi_N^*h(\tilde X,\tilde Y))
		-\pi_N^*h(\widetilde{[V,X]},\tilde Y) - \pi_N^*h(\tilde X,\widetilde{[V,Y]})
	\end{equation*}
	From \Cref{lem:twisted-structures}, we have
	\begin{equation*}
		\widetilde{[V,Y]}=\widehat{[V',Y']}-\pi_N^*(f_Z^{-1}F(V,Y))\cdot \tilde Z
	\end{equation*}
	and therefore we find
	\begin{align*}
		\pi_N^*\big((L_V h)(X,Y)\big)
		&=\big(L_{\widehat{V'}}(\pi_{\bar N}^*h')\big)(\widehat{X'},\widehat{Y'})\\
		&\quad +\pi_N^*(f_Z^{-1}\omega(V,X))\pi_{\bar N}^*h'(\tilde Z,\tilde Y)
		+\pi_N^*(f_Z^{-1}\omega(V,Y))\pi_{\bar N}^*h'(\tilde X,\tilde Z)\\
		&=\pi^*_{\bar N}\big((L_{V'}h')(X',Y')\big)
		+\pi_{N}^*\big(2f_Z^{-1}(\iota_V \omega)\vee (\iota_Z h)\big)(\tilde X,\tilde Y)
	\end{align*}
	Bringing the second term to the left-hand side yields the claimed result.
\end{myproof}

\subsubsection{Elementary deformations}

Given a hyper-K\"ahler manifold with circle action, it is natural to expect that---after choosing suitable twist data---its twist manifold also admits interesting (quaternionic) structures. In this setting, Macia and Swann~\cite{MS2015} investigated when the twist construction yields a quaternionic K\"ahler manifold. They discovered that, in order to obtain a quaternionic K\"ahler metric on the twist manifold, the original hyper-K\"ahler structure must first be modified. To study these modifications, they introduced the notion of elementary deformations of hyper-K\"ahler metrics, which we now discuss.

We start from a (pseudo-)hyper-K\"ahler manifold $(N,g,I_i)$ with K\"ahler forms $\omega_i$ ($i=1,2,3$). For notational convenience, we set $\omega_0=g$ and $I_0=\id$. To apply the twist construction, we need a circle action generated by a vector field $Z$. We will require that the action generated by $Z$ is compatible with the quaternionic Hermitian structure in the sense that $Z$ preserves $g$, as well as the subbundle of  $\End(TM)$ spanned by $I_1$, $I_2$ and $I_3$. We define one-forms $\alpha_i(\cdot )=\omega_i(Z,\cdot )$, $i=0,1,2,3$, as well as the symmetric $(0,2)$-tensor $g_\alpha=\sum_{i=0}^3 \alpha_i^2$. When $Z$ is nowhere-vanishing, this is proportional to the restriction of the metric $g$ to the quaternionic span $\H Z=\langle Z,I_1Z,I_2Z,I_3Z\rangle$ of $Z$.

\begin{mydef}
	Given a pseudo-hyper-K\"ahler manifold $(M,g,I_i)$, an \emph{elementary deformation} of $g$ is a (pseudo-)Riemannian metric on $(M,I_i)$ of the form $g_\h=a\cdot g+ b\cdot g_\alpha$ for nowhere-vanishing $a,b\in C^\infty(M)$.
\end{mydef}

The action of $L_Z$ on $\langle I_1,I_2,I_3\rangle\subset \End TM$ is linear, because $\nabla I_i=0$. Furthermore, the Leibniz rule implies that $L_Z\in \mf{so}(3)$. Therefore, the action of $Z$ is either trivial, or conjugate to
\begin{equation*}
	L_ZI_1=0 \qquad \qquad L_Z I_2=I_3 \qquad \qquad L_Z I_3 =-I_2
\end{equation*}
after suitable normalization of $Z$. In the latter case, we say that $Z$ generates a rotating circle symmetry of $(M,g,I_i)$. In the following, we will consider only such rotating symmetries. In fact, we will make one further assumption, namely that the action generated by $Z$ is $\omega_1$-Hamiltonian: $\iota_Z\omega_1=-\d f_Z$ for some moment map $f_Z \in C^\infty(M)$. Such a hyper-K\"ahler manifold, equipped with $\omega_1$-Hamiltonian rotating circle symmetry, is precisely the input of the HK/QK correspondence.

\subsubsection{The HK/QK correspondence}

It was shown by Haydys~\cite{Hay2008} that, given a rotating circle symmetry on a hyper-K\"ahler manifold, one can construct a bundle over it which admits the structure of a hyper-K\"ahler cone. The hyper-K\"ahler cone then gives rise to a quaternionic K\"ahler manifold by the inverse of the Swann bundle construction~\cite{Swa1991}. His construction is known as the HK/QK correspondence.

These results (which Haydys obtained under the assumption that all metrics involved are positive definite), were extended to metrics of arbitrary signature in \cite{ACM2013}, with control over the signature of the resulting quaternionic K\"ahler metrics. In particular, necessary and sufficient conditions for the resulting metric to be positive definite were given, and shown to include not only definite but also indefinite initial data. In this way, the authors of \cite{ACM2013,ACDM2015} were able to obtain the supergravity $c$-map as a special case of the HK/QK correspondence with indefinite initial data. As part of this work, they gave a simple expression (cf.~\cite[Theorem 2]{ACDM2015}) for the quaternionic K\"ahler metric constructed in the HK/QK correspondence, which is precisely the metric given in the theorem below.

Macia and Swann \cite{MS2015} then showed that the above-mentioned metric can be obtained from a combination of an elementary deformation and a twist, applied to the initial hyper-K\"ahler manifold. Furthermore, they characterized those elementary deformations and twist data that lead to a quaternionic K\"ahler metric. In summary, one obtains the following result:

\begin{thm}[\cite{Hay2008}; \cite{ACM2013}; \cite{ACDM2015}; \cite{MS2015}]\label{thm:deformationcmap}
	Let $(M,g,I_i)$ be a (possibly indefinite) hyper-K\"ahler manifold equipped with rotating circle symmetry $Z$, such that $\omega_1$ is integral and the lift $Z_P$ generates a free and proper action. Denoting the Hamiltonian function of $Z$ with respect to $\omega_1$ by $f_Z$, one can associate a (possibly indefinite) quaternionic K\"ahler manifold $(\bar N,g_\q)$ as follows. Let $g_\h$ be an elementary deformation of $g$ with respect to $Z$, and consider the twist manifold $(\bar N,g_\q)$, where $g_\q\sim_{\mc H} g_\h$ with respect to the twist data $(Z,\omega_\h,f_\h)$. Then $(\bar N,g_\q)$ is quaternionic K\"ahler if and only if 
	\begin{gather*}
		\omega_\h=k(\omega_1+\d \alpha_0)\hspace{1.5cm} f_\h=k(f_Z+g(Z,Z))\\
		g_\h=\frac{B}{f_Z}g+\frac{B}{f_Z^2}g_\alpha 
	\end{gather*}
	where $B,k\in \R\setminus\{0\}$.
\end{thm}

\begin{rem}\leavevmode\label{rem:conventions}
	\begin{numberedlist}
		\item $B$ and $k$ simply scale the metric and curvature of $P$, respectively. In the following, we will set them to $1$. Since the Hamiltonian function $f_Z$ is unique only up to an additive constant, there is an implicit dependence of the construction on this additional parameter. In the case of the $c$-map, this parameter governs the so-called one-loop deformation of the quaternionic K\"ahler metric.
		\item For the convenience of the reader, we provide a dictionary to translate between the notation used for the various objects living on the hyper-K\"ahler manifold in this article, in the work by Macia and Swann \cite{MS2015}, and in the paper \cite{ACDM2015} by Alekseevsky, Cort\'es, Dyckmanns and Mohaupt:
		\begin{center}
			\begin{tabular}{@{}lll@{}} \toprule 
				This article & Macia \& Swann & ACDM \\\midrule 
				$Z$ & $X$ & $-\frac{1}{2}Z$ \\
				$f_Z$ & $-(\mu-c)$ & $-\frac{1}{2}f$\\
				$\eta$ & $\theta$ & $\eta$\\
				$\alpha_0$ & $\alpha_0$ & $-\frac{1}{2}\beta$ \\
				$\omega_\h$ & $F$ & $\omega_1-\frac{1}{2}\d\beta$\\
				$f_\h$ & $-a$ & $\frac{1}{2}f_1$\\
				$g_\h$ & $-g_\h$ & $2g'$\\\bottomrule
			\end{tabular}
		\end{center}
	\end{numberedlist}
\end{rem}

\subsection{Special Kähler geometry and the \texorpdfstring{$c$}{c}-map}
\label{subsec:specialK}

\subsubsection{Special Kähler structures}

In this section, we review the $c$-map and the various manifolds involved in this construction.

\begin{mydef}
	An \emph{affine special K\"ahler manifold} $(M,g,\omega,\nabla)$ is a pseudo-K\"ahler manifold $(M,g,\omega)$ endowed with a flat, torsion-free connection $\nabla$ such that $\nabla \omega=0$ and $\d_\nabla J=0$, where $J$ is the complex structure on $M$.
\end{mydef}

For the $c$-map, the relevant class of affine special K\"ahler manifolds is the following:

\begin{mydef}
	An affine special K\"ahler manifold $(M,g,\omega,\nabla)$ is called \emph{conical} (or a \emph{CASK manifold}) if it admits a vector field $\xi$, called the \emph{Euler field}, which satisfies $\nabla \xi =\nabla^\text{LC}\xi=\id_{TM}$, where $\nabla^\text{LC}$ denotes the Levi-Civita connection of $g$. We will furthermore assume throughout this paper that $\{\xi,J\xi\}$ generates a principal $\C^*$-action, and that $\mc D\coloneqq \Span\{\xi,J\xi\}$ is negative definite with respect to $g$, while the orthogonal complement $\mc D^\perp$ is positive definite.
\end{mydef}

The Euler field is homothetic and holomorphic, while $J\xi$ is holomorphic and Killing (cf.\ Proposition 3 of \cite{CM2009}). 

\begin{lem}\label{lem:Euleraffine}
	The Euler field is always $\nabla$-affine, while $J\xi$ is $\nabla$-affine if and only if $\nabla=\nabla^\text{LC}$.
\end{lem}
\begin{myproof}
	Since $\nabla$ is flat we may use local, parallel vector fields $X,Y$ and compute 
	\begin{equation*}
		(L_\xi \nabla)_XY=L_\xi(\nabla_XY)-\nabla_{[\xi,X]}Y-\nabla_X([\xi,Y])=\nabla_XY=0
	\end{equation*}
	where we used $[\xi,Y]=-\nabla_Y \xi=-Y$. Analogously, we obtain 
	\begin{align*}
		(L_{J\xi}\nabla)_X Y&=-\nabla_X([J\xi,Y])=\nabla_X(\nabla_Y(J\xi))
		=\nabla_X((\nabla_YJ)\xi+J\nabla_Y\xi)\\
		&=\nabla_X((\nabla_\xi J)Y+J\nabla_Y\xi)
		=\nabla_X(JY)=(\nabla_X J)Y
	\end{align*}
	where we used $\nabla_\xi J=(L_\xi+\nabla \xi)J=0$ and $\nabla\xi=\id_{TM}$ in the penultimate step. This expression vanishes if and only if $\nabla$ is the Levi-Civita connection.
\end{myproof}

The following fact is well-known (see, for instance, \cite{CHM2012,MS2015}).

\begin{lem}\label{lem:CASKpotential}
	The function $f=\frac{1}{2}g(\xi,\xi)$ defines both a Hamiltonian function for $J\xi$ and a K\"ahler potential for $g_M$.
\end{lem}
\begin{myproof}
	$\nabla \xi=\id_{TM}$ implies $\d f=g(\xi,\cdot)=-\omega(J\xi,\cdot)$. Furthermore $\d \d^cf=-\d J^*\d f=\d(\omega (J\xi,J\cdot))=\d(\omega(\xi,\cdot))=2\omega$, where we used $\nabla\omega=0$ and $\nabla \xi=\id_{TM}$.
\end{myproof}

Since $\xi$ generates a homothetic $\R_{>0}$-action, dividing out the $\C^*$-action on a CASK manifold amounts to fixing a level set for the moment map $\frac{1}{2}g(\xi,\xi)$---we will take the level where $g(\xi,\xi)=-1$ for notational convenience---and subsequently taking the quotient by the circle action generated by $-J\xi$.

\begin{mydef}
	A \emph{projective special K\"ahler} (PSK) manifold $\bar M$ is the K\"ahler quotient $M\sslash S^1$ of a CASK manifold $M$ by the $S^1$-action generated by the Hamiltonian Killing field $-J\xi$.
\end{mydef}

The natural notion of symmetry in this setting is the following:

\begin{mydef}\leavevmode
	\begin{numberedlist}
		\item An automorphism of a CASK manifold $(M,g_M,\omega_M,\nabla,\xi)$ is a diffeomorphism of $M$ which preserves $g_M$, $\omega_M$, $\nabla$ and $\xi$.
		\item Let $\bar M=M\sslash S^1$ be a PSK manifold. An automorphism of $\bar M$ is a diffeomorphism of $\bar M$ induced by an automorphism of the CASK manifold $M$.
	\end{numberedlist}
\end{mydef}

\begin{lem}\label{lem:CASKHam}
	If $X\in \mf X(M)$ generates a one-parameter group $\varphi_t$ of automorphisms of a CASK manifold $(M,g_M,\omega_M,\nabla,\xi)$, then $X$ is $\omega_M$-Hamiltonian.
\end{lem}
\begin{myproof}
	Defining $\alpha\coloneqq \frac{1}{2}\d^c(g(\xi,\xi))\in \Omega^1(M)$, \Cref{lem:CASKpotential} asserts $\omega_M=\d\alpha$. Then $\rho\coloneqq \alpha(X)$ is a candidate Hamiltonian function. We compute its differential: 
	\begin{equation*}
		\d \rho=\d\iota_X\alpha=(L_X-\iota_X\d)\alpha=L_X\alpha-\iota_X\omega_M
	\end{equation*}
	Thus, it suffices to check that $L_X\alpha=0$. Since $\alpha=-\frac{1}{2}J^*(\d g(\xi,\xi))$, this follows from the fact that $X$ preserves $J$, $g$ and $\xi$.
\end{myproof}

\begin{mydef}\leavevmode
	\begin{numberedlist}
		\item An infinitesimal automorphism of a CASK manifold $(M,g_M,\omega_M,\nabla,\xi)$ is a vector field $X\in M$ such that its local flow preserves the CASK data on $M$. The Lie algebra of such vector fields is denoted by $\aut(M)$.
		\item An infinitesimal automorphism of a PSK manifold is a vector field $\bar X$ induced by an infinitesimal automorphism of the corresponding CASK manifold (note that the latter always project, since they commute with $\xi$ and $J\xi$). The corresponding Lie algebra is denoted by $\aut(\bar M)$.
	\end{numberedlist}
\end{mydef}

\begin{prop}\label{lem:infPSKCASKaut}
	Consider the natural map $\varphi:\aut(M)\to \aut(\bar M)$ induced by the projection $\pi:M\to \bar M$.
	\begin{numberedlist}
		\item If $\nabla\neq \nabla^\text{LC}$, then $\varphi$ is an isomorphism, with inverse $\bar X\mapsto X^{\mc H}-\pi^*h_{\bar X}J\xi$, where $X^{\mc H}$ denotes the horizontal (that is, perpendicular to $\mc D$) lift of $\bar X$ and $h_{\bar X}$ is a uniquely determined Hamiltonian for $\bar X$ with respect to the K\"ahler form $\omega_{\bar M}$ on $\bar M$.
		\item If $\nabla=\nabla^\text{LC}$, then $\varphi$ induces an isomorphism between $\aut(M)/(\R \cdot  J\xi)$ and $\aut(\bar M)$, with inverse $\bar X\mapsto X^{\mc H} -\pi^*h_{\bar X}J\xi$, where $h_{\bar X}$ is now any Hamiltonian for $\bar X$ with respect to $\omega_{\bar M}$
	\end{numberedlist}
\end{prop}
\begin{myproof}
	Since the $\C^*$-action generated by $\{\xi,J\xi\}$ is homothetic, the distribution $\mc D^\perp$ of vectors orthogonal to the vertical distribution $\mc D=\Span\{\xi,J\xi\}$ defines a principal connection in the $\C^*$-bundle $M\to\bar M$.
	
	Now consider an infinitesimal automorphism $X\in \aut(M)$, which projects to $\bar X\in \aut(\bar M)$. We can decompose $X=X^{\mc H}+X^{\mc V}$ into its horizontal and vertical components, where $X^{\mc H}$ is the horizontal lift of $\bar X$. Since $X$ preserves both $g_M$ and $\xi$, we have $0=L_X(g_M(\xi,\xi))=2g_M(\nabla^\text{LC}_X\xi,\xi)=2g_M(X,\xi)$, i.e.\ $X$ is orthogonal to $\xi$.
	
	Thus, we may write $X^{\mc V}=-f_XJ\xi$ for some function $f_X\in C^\infty(M)$. Since $X^{\mc H}$ is $\C^*$-invariant, we obtain that $0=[X,\xi]=[X^{\mc V},\xi]=\xi(f_X)J\xi$ and similarly for $J\xi$. This shows that $f_X$ is constant on the fibers of the projection, hence of the form $\pi^*h_{\bar X}$ for some $h_{\bar X}\in C^\infty(\bar M)$. We now show that $h_{\bar X}$ is uniquely determined if $\nabla\neq \nabla^\text{LC}$, and up to additive constant in case $\nabla=\nabla^\text{LC}$. Indeed, let $X$ and $X'$ be two infinitesimal CASK automorphisms projecting to $\bar X$. Then $X'-X=\pi^*(h_{\bar X}-h'_{\bar X})J\xi$ must be both Killing and $\nabla$-affine. The first condition implies that $h_{\bar X}-h'_{\bar X}$ is constant. In case $\nabla\neq \nabla^\text{LC}$, the latter then implies that it vanishes (by \Cref{lem:Euleraffine}).
	
	We claim that $h_{\bar X}$ is a Hamiltonian of $\bar X$ with respect to the K\"ahler form $\omega_{\bar M}$ of $\bar M$. The principal connection $\mc D^\perp$ induces a principal $S^1$-connection on the bundle $S=\{g(\xi,\xi)=-1\}\subset M$ over $\bar M$, which has $-J\xi$ as its fundamental vector field. We denote its connection one-form by $\theta=g(J\xi,\cdot)$. Differentiating this equation, we obtain $\d \theta=\omega_M|_S=\pi_S^*\omega_{\bar M}$, where $\pi_S:S\to \bar M$ is the projection map. Since $X$ is an infinitesimal CASK automorphism, $L_X\theta=0$, which means that
	\begin{equation*}
		0=\d\iota_X\theta+\iota_X\d\theta=\pi^*\d h_{\bar X} +(\pi^*\omega_M)(X,\cdot)=\pi^*(\d h_{\bar X}+\omega_{\bar M}(\bar X,\cdot))
	\end{equation*}
	This finishes the proof.
\end{myproof}

An extrinsic perspective on special K\"ahler geometry is developed in \cite{ACD2002}. It is proven that any simply connected affine special K\"ahler manifold $M$ of complex dimension $n+1$ admits a (holomorphic, non-degenerate, Lagrangian) immersion into $T^*\C^{n+1}$ such that the standard special K\"ahler structure on $T^*\C^{n+1}$ induces the given affine special K\"ahler structure on $M$. In fact, such an immersion is given by the graph of an exact one-form. This implies that the special K\"ahler structure is determined by a holomorphic function on an open subset of $\C^{n+1}$ satisfying a certain non-degeneracy condition. Thus, any simply connected special K\"ahler manifold may be regarded as a domain $U\subset \C^{n+1}$ endowed with a holomorphic function $F$ which determines the special K\"ahler structure completely. $F$ is known as the holomorphic prepotential of $U$. An analogous theorem holds for CASK manifolds, requiring that the domain and the holomorphic function are conical in the appropriate sense, and by projectivizing one obtains an identification of the corresponding PSK manifold $\bar M$ with an open set in $\CP^n$, inducing the PSK structure. 

More generally, any affine special K\"ahler manifold $M$ can be covered by open sets $U_\alpha$, called special K\"ahler domains, on which the affine special K\"ahler structure arises in this fashion. By picking suitable global coordinates $\{z_I\}_{I=0,1,\dots,n}$ on $\C^{n+1}$ and possibly shrinking the open sets $U_\alpha$, we may assume that (the image of) each $U_\alpha$ is contained in $\{z_0\neq 0\}\subset\C^{n+1}$. From now on, we will only consider special K\"ahler domains of this type. 

In case $M$ is a CASK manifold, projectivizing CASK domains yields local identifications of the corresponding PSK manifold $\bar M$ with subsets of $\{z_0\neq 0\}\subset \CP^n$. To each such PSK domain $\bar U_\alpha$ we may canonically associate a K\"ahler potential, defined analogously to the standard potential for the Fubini--Study metric on $\CP^n$. Recall from \Cref{lem:CASKpotential} that $\frac{1}{2}g(\xi,\xi)$ defines a global K\"ahler potential on $M$. Identifying the CASK and PSK domains with their imagines in $\C^{n+1}$ and $\CP^n$, the function
\begin{equation*}
	\mc K_\alpha\coloneqq \log\bigg(\frac{\frac{1}{2}g(\xi,\xi)}{\abs{z_0}^2}\bigg)
\end{equation*}
now defines a K\"ahler potential on $\bar U_\alpha$.

\begin{rem}
	We revisit the case $\nabla=\nabla^\text{LC}$ from the previous lemma. Then, the short exact sequence
	\begin{equation*}
	\begin{tikzcd}
	0 \ar[r] & \R\cdot J\xi \ar[r] & \aut(M) \ar[r] & \aut(\bar M) \ar[r] & 0
	\end{tikzcd}
	\end{equation*}
	of Lie algebras splits. To see this, let us first assume that $M$ is simply connected. Then, in light of the above discussion, we can (globally) realize it as a complex Lagrangian submanifold of $T^*\C^{n+1}$ (where $n=\dim_\C \bar M$). From the condition $\nabla=\nabla^\text{LC}$, it follows that it is contained in a linear subspace, and that there exists a global, holomorphic prepotential of the form $F=\sum a_{ij}z^iz^j$, where the real matrix $(N_{ij})\coloneqq 2\Im (a_{ij})$ is of signature $(n,1)$. In fact, the condition $\nabla=\nabla^\text{LC}$ is equivalent to the vanishing of the third derivatives of $F$ \cite{Fre1999}. This implies that $\aut (M)=\mf{u}(n,1)$ and $\aut (\bar M)=\mf{su}(n,1)$, with $\R\cdot J\xi=\mf{u}(1)$ the center of $\mf{u}(n,1)$, so the sequence splits. In case $\pi_1(M)=\Gamma$ is non-trivial, we can write $M=\tilde M/\Gamma$, where $\tilde M$ is the universal covering, which is also a CASK manifold. Then we have $\R\cdot J\xi=\mf{u}(1)\subset \aut(M)=\aut(\tilde M)^\Gamma\subset \aut(\tilde M)=\mf{u}(n,1)=\mf{u}(1)\oplus \mf{su}(n,1)$, which implies that $\aut(M)=\mf{u}(1)\oplus \mf{h}$, where $\mf{h}\subset \mf{su}(n,1)$.
\end{rem}

\subsubsection{The rigid and supergravity \texorpdfstring{$c$}{c}-maps}

As mentioned in the introduction, the (supergravity) $c$-map and its one-loop deformation construct a complete quaternionic K\"ahler manifold $\bar N$ of negative scalar curvature out of any complete PSK manifold $\bar M$. We first describe it under the assumption that $\bar M$ is a PSK domain \cite{CDS2017}. Then, we may regard $\bar M$ as a subset of $\CP^n$ endowed with a canonical K\"ahler potential $\mc K$, and the corresponding CASK manifold $M$ as a domain in $\C^{n+1}$ equipped with a holomorphic function $F$.

As a smooth manifold, $\bar N$ is simply the direct product $\bar M\times \R^{2n+4}$, but its Riemannian structure given by the one-loop deformed Ferrara--Sabharwal metric $g_\text{FS}^c$ is not that of a product. To write this metric down explicitly, we use global coordinates. On $\bar M\subset \CP^n$ we use homogeneous coordinates $(z^0:z^1:\dots:z^n)$, which under the standard identification $\{z^0\neq 0\}\cong \C^n$ correspond to coordinates $(X^0\equiv 1,X^1,\dots,X^n)$ on $\C^n$. Though $X^0\equiv 1$, it will nevertheless appear in the following expressions to simplify notation. We further have coordinates $(\rho,\tilde\phi,\zeta^I,\tilde \zeta_I)$ on $\R_{>0}\times \R\times \R^{n+1}\times \R^{n+1}\cong \R^{2n+4}$. With respect to these, and denoting the K\"ahler metric on $\bar M$ by $g_{\bar M}$, $g_\text{FS}^c$, $c\geq 0$, is given by
\begin{equation}\label{eq:gFSstructure}
	g_\text{FS}^c=\frac{\rho+c}{\rho}g_{\bar{M}}+g^c_G
\end{equation}
where 
\begin{align*}
	g_G^c&\coloneqq \frac{1}{4\rho^2}\frac{\rho+2c}{\rho+c}\d\rho^2
	+\frac{1}{4\rho^2}\frac{\rho+c}{\rho+2c}\Big(\d\tilde\phi+\sum_{I=0}^n(\zeta^I\d \tilde \zeta_I-\tilde\zeta_I\d\zeta^I)+c\d^c\mc K\Big)^2\\\numberthis\label{eq:FSmetric}
	&+\frac{1}{2\rho}\sum_{I,J=0}^n\bigg(\mc J_{IJ} \d \zeta^I \d \zeta^J
	+\mc J^{IJ}(\d \tilde \zeta_I+\mc R_{IK}\d \zeta^K)(\d \tilde \zeta_J + \mc R_{JL}\d \zeta^L)\bigg)\\
	&+\frac{2c}{\rho^2}e^{\mc K}\Big|\sum_{I=0}^n(X^I\d \tilde\zeta_I+F_I\d \zeta^I)\Big|^2
\end{align*}
In this expression, $(F_I)$, $(\mc J_{IJ})$ and $(\mc R_{IJ})$ are determined by the derivatives of the holomorphic prepotential $F$, and constant on each copy of $\R^{2n+4}$. Their precise definitions can be found in \cite{CDS2017}, but will not be relevant in the following. We note that, for any two $c_1,c_2>0$, the metrics $g^{c_1}_\text{FS}$ and $g^{c_2}_\text{FS}$ are isometric~\cite{CDS2017}.

It was explained in \cite{CDS2017} how to define the $c$-map for arbitrary PSK manifolds. In the case $c=0$ the construction is locally identical to the $c$-map for special K\"ahler domains, with the only global difference being that $\bar N$ may fail to be globally trivial as a bundle over $\bar M$. This corresponds to the fact that, after covering $\bar M$ by PSK domains, transition functions must be used to patch the different domains together. 

For $c>0$, however, this patching does not immediately preserve the metric $g^c_\text{FS}$: One must first divide out the isometric $\Z$-action which sends $\tilde\phi\mapsto \tilde\phi+2\pi c$, thereby making the coordinate $\tilde \phi$ periodic. Thus, in order to define the one-loop deformed $c$-map on arbitrary PSK manifolds, one locally works with bundles with fiber $\R_{>0}\times S^1\times \R^{2n+2}$ instead. Then, the metric $g^c_\text{FS}$ is invariant under the transition functions, and therefore patches together to a global quaternionic K\"ahler metric. Note that its coordinate expression does not change.

The resulting manifold $\bar N$ is best thought of as a bundle of Lie groups over the PSK manifold $\bar M$. If $\bar M$ is a PSK domain, we may endow each fiber with a group structure by casting $\R^{2n+4}\cong \R_{>0}\times \R^{2n+3}$ as the Iwasawa subgroup of $SU(n+2,1)$, which we will denote by $G(n+2)$. This is a one-dimensional, solvable extension of the Heisenberg group of dimension $2n+3$, which is parametrized by the coordinates $(\tilde\phi,\tilde\zeta_I,\zeta^I)$ in the above notation. With respect to these coordinates, the group multiplication takes the following form:
\begin{align*}
	(e^\lambda,\alpha,\tilde v_I,v^I)\cdot &(\rho,\tilde\phi,\tilde \zeta_I,\zeta^I)\\\numberthis \label{eq:Gaction}
	&\quad =\Big(e^\lambda \rho,
	e^{\lambda}\tilde\phi+\alpha+e^{\lambda/2}\Big(\sum_I \tilde v_I\zeta^I-v^I\tilde \zeta_I\Big),
	\tilde v_I+e^{\lambda/2}\tilde \zeta_I,v^I+e^{\lambda/2}\zeta^I\Big)
\end{align*}
In the case of PSK domains (or if $c=0$), this action endows $\bar N$ with the structure of a principal $G(n+2)$-bundle. In the general case (for $c>0$), when the $\tilde\phi$-coordinate is periodic, one no longer has an action of all of $G(n+2)$ on the fibers, since the one-dimensional extension fails to respect the quotient. However, shifts in $\tilde \phi$ are central in $\Heis_{2n+3}$, hence its action descends to $\bar N$, and we may take the appropriate cyclic quotient to obtain a free action. Regarding the compatibility of the Ferrara--Sabharwal metric and its one-loop deformation with these actions, we have the following basic result.

\begin{lem}\label{lem:fiberisometries}
	Consider a $c$-map space $\bar N$ equipped with the the (one-loop deformed) Ferrara--Sabharwal metric $g^c_\text{FS}$. Then $\Heis_{2n+3}$ acts by isometries for every $c\geq 0$, and all of $G(n+2)$ acts by isometries in case $c=0$.
\end{lem}
\begin{myproof}
	Let $U\subset \bar M$ be a PSK domain, so that $\bar N$ restricts to a trivial bundle over $U$. Using the notation introduced above, the metric $g^c_\text{FS}$ then takes the form given by equations \eqref{eq:gFSstructure} and \eqref{eq:FSmetric}, and the action of $G(n+2)$ is given by \eqref{eq:Gaction}. In case $\bar M$ was not globally a PSK domain and $c>0$, only the subgroup $\Heis_{2n+3}$ has a well-defined action on $\bar N$. Since the group acts purely in the fiber directions, it is enough to check that its elements leave $g^c_G$ invariant and preserve the function $\frac{\rho+c}{\rho}$ to verify that they act by isometries. That this is the case for elements of the form $(0,\alpha,\tilde v_I,v^I)\in \Heis_{2n+3}\subset G(n+2)$ follows easily from an explicit computation. 
	
	The remaining elements, which are of the form $(e^\lambda,0,0,0)\in G(n+2)$, only act on $\bar N$ when either $c=0$ or $\bar M$ is globally a PSK domain (and $c>0$). In the former case, they act by isometries. In the latter, we observe that for $\lambda\neq 0$, these elements preserve neither the function $\frac{\rho+c}{\rho}$ nor $g^c_G$, and therefore do not induce isometries.
\end{myproof}

\begin{cor}
	Let $\bar M$ be a complete, simply connected PSK manifold of (real) dimension $2n$, and $\Gamma\subset \Heis_{2n+3}$ a cocompact lattice. Then the (universal covering of the) one-loop deformed $c$-map space $(\bar N,g^c_\text{FS})$, $c\geq 0$, associated to $\bar M$ admits a smooth quotient $\bar N/\Gamma$, which is a complete quaternionic K\"ahler manifold with fundamental group $\Gamma$.
\end{cor}
\begin{myproof}
	If either $c=0$ or $\bar M$ is a PSK domain, $\bar N$ is itself simply connected, and carries a free and isometric action of $\Heis_{2n+3}$. In case $\bar M$ is not a PSK domain and $c>0$, this is still true for its universal covering (which we also denote by $\bar N$ for notational convenience). The quotient by $\Gamma$ is smooth, complete, and has fundamental group $\Gamma$ by construction.
\end{myproof}

This direct description gives a satisfactory (global) definition of the $c$-map, but makes it difficult to study the properties of the $c$-map. For instance, there is no obvious way to compute the isometry group of $(\bar N,g^c_\text{FS})$ other than by direct inspection of the coordinate expressions. An alternative, locally equivalent, formulation is summarized by \eqref{eq:diagram}, which we now explain in more detail. 

The cotangent bundle of an affine special K\"ahler manifold carries a (pseudo-)hyper-K\"ahler structure which admits a rather simple description. Let $(M,g_M,\omega_M,\nabla)$ be a special K\"ahler manifold. The special K\"ahler connection induces a splitting 
\begin{equation}\label{eq:splitting}
	T(T^*M)\cong T^{\mc H}N\oplus T^{\mc V}N\cong \pi^*(TM)\oplus \pi^*(T^*M)
\end{equation} 
With respect to this splitting, we define
\begin{equation}\label{eq:cotangentHKstr}
	g=
	\begin{pmatrix}
		g_M & 0 \\ 0 & g^{-1}_M
	\end{pmatrix}
	\qquad  
	I_1=
	\begin{pmatrix}
		J_M & 0 \\ 0 & J_M^*
	\end{pmatrix}
	\qquad  
	I_2=
	\begin{pmatrix}
		0 & -\omega_M^{-1} \\ \omega_M & 0
	\end{pmatrix}
	\qquad 
	I_3=I_1I_2
\end{equation}
where we omitted pullbacks to simplify the notation. A priori, these tensor fields (only) endow $N\coloneqq T^*M$ with an almost hyper-Hermitian structure. However, $(N,g,\omega_a)$ is in fact hyper-K\"ahler if (and only if) $(M,g_M,\omega_M,\nabla)$ is special K\"ahler (see, for instance, \cite{ACD2002}). In case $M$ is a CASK manifold, the map which sends $M$ to $N$, equipped with the above hyper-K\"ahler structure, is called the rigid $c$-map. It is also well-known that, in this case, $N$ carries a natural rotating circle symmetry (cf.\ \cite{ACM2013}):

\begin{lem}
	Let $(M,g_M,\omega_M,\nabla,\xi)$ be a CASK manifold, and $(N,g,I_i)$, $i=1,2,3$, its cotangent bundle, equipped with the rigid $c$-map hyper-K\"ahler structure. Then the vector field $Z=-\widetilde{J\xi}$ on $N$, where the tilde denotes the $\nabla$-horizontal lift, generates a rotating circle symmetry. In other words, $L_Zg=L_Z\omega_1=0$, $L_Z\omega_2=\omega_3$ and $L_Z\omega_3=-\omega_2$. \proofclear 
\end{lem}

Since $Z$ is horizontal, $\omega_1(Z,\cdot )=\pi^*\omega_M(Z,\cdot )$. By \Cref{lem:CASKpotential}, the rotating circle symmetry is Hamiltonian with Hamiltonian function $-\frac{1}{2}\pi^*(g_M(\xi,\xi))=-\frac{1}{2}g(Z,Z)$. As explained in the previous section, one may now apply the HK/QK correspondence to obtain the Ferrara--Sabharwal metric~\cite{ACDM2015,MS2015}.\footnote{Due to global issues when dividing out the lifted circle action on the $S^1$-bundle over $N$, the HK/QK correspondence may produce a manifold which is not globally diffeomorphic to the manifold $\bar N$ obtained by the direct $c$-map construction given above. Nevertheless, the approaches are locally equivalent (as can be seen by picking a local transversal slice for the action generated by $Z_P$), and therefore we may use the HK/QK correspondence to study tensor fields on $\bar N$.} Though $-\frac{1}{2}g(Z,Z)$ is a natural choice of Hamiltonian function for $Z$ with respect to $\omega_1$, it is not unique. We may choose to work with the function $f_Z=-\frac{1}{2}g(Z,Z)-\frac{1}{2}c$ for any $c\in\R$ instead (the normalization of $c$ is conventional). This induces a change in the twist data used in the HK/QK correspondence, since $f_\h=f_Z+g(Z,Z)$ is shifted as well. Consequently, there is a one-parameter freedom in the HK/QK correspondence (cf.\ \Cref{rem:conventions}), which precisely reproduces the one-loop deformation of the Ferrara--Sabharwal metric~\cite{MS2015}.

There are two main advantages to the HK/QK correspondence over the direct description of the $c$-map. Firstly, since the one-loop deformation of the $c$-map corresponds to a simple shift in the Hamiltonian function $f_Z$ on the hyper-K\"ahler side of the correspondence, it is much simpler to treat the undeformed $c$-map simultaneously with its one-loop deformation. Secondly, the twist construction makes it possible to compute tensor fields on the $c$-map image $\bar N$ in terms of their counterparts on the hyper-K\"ahler manifold $N$. We will exploit both of these in the following.

\clearpage

\section{Symmetry properties of the HK/QK correspondence and the \texorpdfstring{$c$}{c}-map}
\label{sec:symmetries}

We now investigate the behavior of the HK/QK correspondence and the $c$-map in the presence of symmetries. We will show that under the supergravity $c$-map as well as under its one-loop deformation, the identity (connected) component of the automorphism group of the initial projective special K\"ahler manifold is a subgroup of the isometry group of the resulting quaternionic K\"ahler manifold. In the undeformed case, this result was already known (see e.g.\ \cite[App.~A]{CDJL2017}), but our usage of the HK/QK correspondence ensures that our results remain valid also for the one-loop corrected supergravity $c$-map.

In addition to isometries of the quaternionic K\"ahler manifold coming from the automorphisms of the PSK manifold, there is a solvable group $G(n+2)$ of dimension $2n+4$ acting on the quaternionic K\"ahler manifold, with orbits transverse to the orbits of the automorphism group of the PSK manifold (see, for example, \cite{CHM2012}). In the undeformed case, it acts isometrically. In the deformed case, however, only its nilradical, a codimension one Heisenberg subgroup, acts by isometries (see \Cref{thm:dimisom}).

\subsection{Automorphisms under the HK/QK correspondence}
\label{sec:HKQKaut}

As explained in \Cref{sec:prelims}, the input for the HK/QK correspondence is a (connected) hyper-K\"ahler manifold $(N,g,I_i)$ ($i=1,2,3$) equipped with a so-called rotating symmetry, i.e.\ a circle action which is isometric and $\omega_1$-Hamiltonian---with Hamiltonian function $f_Z$---and whose generating vector field $Z$ satisfies $L_Z\omega_2=\omega_3$, $L_Z\omega_3=-\omega_2$.

We now want to consider this set-up in the presence of additional symmetries. Since we are interested in the connected component of the identity of the appropriate symmetry groups, we may work infinitesimally. 

\begin{mydef}
	By an \emph{(infinitesimal) symmetry} of a hyper-K\"ahler manifold with rotating symmetry, we understand a vector field $Y\in \mf X(N)$ which is Killing and triholomorphic as well as $\omega_1$-Hamiltonian, and preserves the Hamiltonian function $f_Z$. The space of such vector fields will be denoted by $\aut(N,f_Z)$.
\end{mydef}

We note that this implies that $[Y,Z]=0$, since $Z=-\omega_1^{-1}(\d f_Z)$, where we interpret $\omega_1$ as an isomorphism $TM\to T^*M$.

Our aim is to produce a Killing field $Y_\q$, associated with $Y$, on the quaternionic K\"ahler manifold $\bar N$, which results from applying the HK/QK correspondence to $N$ with the above data. The first step is to apply the elementary deformation, chosen according to \Cref{thm:deformationcmap}, such that the deformed metric $g_\h=\frac{1}{f_Z}g+\frac{1}{f_Z^2}g_\alpha$ is $\mc H$-related to a quaternionic K\"ahler metric $g_\q$ on the twist manifold $\bar N$. The relevant twist data are $(Z,\omega_\h,f_\h)$, defined as in \Cref{thm:deformationcmap} (we set $k=1$, cf.~\Cref{rem:conventions}).

\begin{lem}
	The vector field $Y$ is $\omega_\h$-Hamiltonian, and preserves the elementary deformation $g_\h$ of the metric. The corresponding Hamiltonian function $f_Y$ (unique up to an additive constant) is $Z$-invariant.
\end{lem}
\begin{myproof}
	We assumed that $Y$ is $\omega_1$-Hamiltonian; call its Hamiltonian function $\varphi$. Then we have
	\begin{equation*}
		\iota_Y\omega_\h=-\d\varphi +\iota_Y\d\alpha_0=-\d (\varphi+g(Z,Y))
	\end{equation*}
	where we used Cartan's formula and the fact that $Y$ preserves both $g$ and $Z$.
	
	It follows directly from our assumptions $L_Yg=0$, $L_YZ=0$, $L_Yf_Z=0$ and $L_Y\omega_i=0$ ($i=1,2,3$) that $Y$ also preserves $g_\h=\frac{1}{f_Z}g+\frac{1}{f_Z^2}g_\alpha$.
	
	To check $Z$-invariance of $f_Y\coloneqq \varphi+g(Z,Y)$, we note $L_Z(g(Z,Y))=0$ because $Z$ is Killing and commutes with $Y$. Thus, it suffices to calculate $L_Z\varphi=\iota_Z \d \varphi=-\omega_1(Y,Z)=-\d f_Z(Y)=0$.
\end{myproof}

Since $Y$ is $Z$-invariant, we can now use the twist construction to define an $\mc H$-related vector field $Y'$ on the (quaternionic K\"ahler) twist manifold $\bar N$. Since we want to produce an infinitesimal isometry of the quaternionic K\"ahler metric $g_\q$, we use \Cref{lem:twisted-Killing} to see that $Y'$ will be Killing with respect to $g_\q$ precisely if
we have 
\begin{equation}\label{eq:Killingtwist}
	L_Y g_\h-2f_\h^{-1}(\iota_Y \omega_\h)\vee (\iota_Z g_\h)=0 
\end{equation}
Since $L_Y g_\h=0$ and the second term does not generally vanish, we will modify $Y$ (in a $Z$-invariant way!) to ensure its twist will be a Killing field. A natural choice is to add a term proportional to the distinguished vector field $Z$. The following shows that this is indeed the right approach:

\begin{lem}
	There exists a $Z$-invariant function $\psi\in C^\infty(N)$ such that the twist of $Y_\h\coloneqq Y+\psi Z$ is Killing.
\end{lem}
\begin{myproof}
	We write $Y_\h=Y+\psi Z$ for some arbitrary $\psi\in C^\infty(N)$, and attempt to solve \Cref{eq:Killingtwist}. Since both $Y$ and $Z$ preserve $g_\h$, we find $L_{Y_\h}g_\h=L_{\psi Z}g_\h
	=2(\d\psi)\iota_Z g_\h$. Since $Y$ and $Z$ are both $\omega_\h$-Hamiltonian, we have:
	\begin{equation*}
		L_{Y_\h}g_\h
		-2f^{-1}_\h(\iota_{Y_\h}\omega_\h)\vee (\iota_Zg_\h)
		=\Big(2(\d\psi)-2f_\h^{-1}(-\d f_Y-\psi \d f_\h)\Big)\vee \iota_Z g_\h
	\end{equation*}
	and thus it suffices, in order for this expression to vanish, to solve the differential equation
	\begin{equation*}
		\d\psi=-\frac{1}{f_\h}\d f_Y-\frac{\psi}{f_\h} \d f_\h
	\end{equation*}
	This equation is solved by $\psi=-\frac{f_Y}{f_\h}$, which therefore satisfies our requirements. The function $\psi$ is $Z$-invariant because both $f_Y$ and $f_\h$ are.
\end{myproof}

Thus, the twist of the vector field $Y_\h=Y-f_\h^{-1}f_YZ$ yields a Killing field $Y'_\h\eqqcolon Y_\q$ of $g_\q$. Recall that $Z_{\bar N}$ is the Killing field on $\bar N$ which plays a role analogous to $Z$ under the HK/QK correspondence.

\begin{lem}\label{lem:YQZNcommute}
	The above-constructed Killing field $Y_\q$ on $\bar N$ commutes with $Z_{\bar N}$.
\end{lem}
\begin{myproof}
	Since $Y_\q\sim_{\mc H} Y_\h$, its horizontal lift to $P$ is $X_P$-invariant---otherwise $Y_\h$ would not be well-defined. The claim now follows from compatibility of the Lie bracket with the pushforward.
\end{myproof}

Denoting the space of Killing fields on $\bar N$ that commute with $Z_{\bar N}$ by $\aut(\bar N,Z_{\bar N})$, we have now proven:

\begin{prop}\label{prop:symmetrytwist}
	Given any $Y\in \mf{aut}(N,f_Z)$, the vector field $Y_\q=\big(Y-\frac{f_Y}{f_\h}Z\big)'$ is an element of $\aut(\bar N,Z_{\bar N})$. \proofclear
\end{prop}

\begin{rem}\label{rem:notunique}
	Our choice of a Hamiltonian function $f_Y$ introduces some non-uniqueness in the definition of $Y_\q$. A different choice $\hat f_Y\coloneqq f_Y+C$ for some $C\in \R$ leads to a new vector field $\hat Y_\q=Y_\q-\frac{C}{f_Z}Z'=Y_\q+CZ_{\bar N}$, where we used that $\tilde Z=Z_P-f_Z X_P$, hence $(\pi_{\bar N})_*\tilde Z=-f_Z Z_{\bar N}$. In summary, $Y_\q$ is unique up to a constant multiple of $Z_{\bar N}$.
\end{rem}

Running through our argument for linearly independent vector fields, we find:

\begin{prop}\label{prop:autmaps}
	Given a basis $\{Y_j\}$, $j=1,\dots,\dim \aut(N,f_Z)$, of $\aut(N,f_Z)$, any set of choices of Hamiltonian functions $f_{Y_j}$ with respect to $\omega_\h$ gives rise to an injective, linear map
	\begin{equation*}
		\begin{tikzcd}[row sep=0cm]
			\varphi:\aut(N,f_Z) \ar[r] & \aut(\bar N,Z_{\bar N}) \\
			Y=\sum_j\alpha_j Y_j \ar[r,mapsto] & Y_\q=\sum_j \alpha_j Y^\q_j
		\end{tikzcd}
	\end{equation*}
	where $Y^\q_j\coloneqq Y_{j,\q}=\Big(Y_j-\frac{f_{Y_j}}{f_\h}Z\Big)'$.
\end{prop}
\begin{myproof}
	Linearity of the map is clear. To prove injectivity, note that the process of twisting a vector field consists of the composition of two point-wise isomorphisms, and therefore introduces no kernel. Thus, $\varphi(Y)=0$ if and only if $Y_\h=0$, i.e.\ $Y=\frac{\sum_j \alpha_j f_{Y_j}}{f_\h}Z$. Any vector field of the form $\psi Z$ for $\psi\in C^\infty(N)$ satisfies $L_{\psi Z}g=2(\d\psi)\iota_Zg$, hence is Killing if and only if $\psi$ is constant. But since $Z$ generates a rotating circle symmetry, no non-zero multiple of it preserves $\omega_2$ or $\omega_3$. Because $Y$ is triholomorphic, $\psi\equiv 0$, so $\ker \varphi=\{0\}$.
\end{myproof}

Both $\aut(N,f_Z)$ and $\aut(\bar N,Z_{\bar N})$ naturally come equipped with the structure of a Lie algebra. It is natural to ask if (any of) the maps from \Cref{prop:autmaps} are Lie algebra homomorphisms. It turns out that this is generally not quite the case:

\begin{thm}\label{thm:structure}
	Let $\{Y_j\}$, $j=1,\dots,k$ be a basis of a Lie subalgebra $\mf g\subset \aut (N,f_Z)$ satisfying $[Y_j,Y_k]=\sum_l c^l_{jk}Y_l$. Then 
	\begin{equation*}
		[Y^\q_j,Y^\q_k]=\sum_lc^l_{jk}Y^\q_l + A_{jk}Z_{\bar N}
	\end{equation*}
	for constants $A_{jk}=\omega_\h'(Y'_j,Y'_k)-\sum_l c^l_{jk} f_{Y_l}'$. Thus, $\mf g$ induces a $(k+1)$-dimensional Lie subalgebra $\mf{g}_\q=\Span \{Y^\q_j,Z_{\bar N}\mid j=1,\dots,k\}\subset\aut(\bar N,Z_{\bar N})$, which is a (possibly trivial) central extension of $\mf g$.
\end{thm}
\begin{myproof}
	Expanding in terms of the $Y_j$, we have:
	\begin{equation*}
		[Y^\h_j,Y^\h_k]=\sum_l c^l_{jk}Y_l
		-\bigg(\bigg[\frac{f_{Y_j}}{f_\h}Z,Y_k\bigg]+\bigg[Y_j,\frac{f_{Y_k}}{f_\h}Z\bigg]\bigg)
		+\bigg[\frac{f_{Y_j}}{f_\h}Z,\frac{f_{Y_k}}{f_\h}Z\bigg]
	\end{equation*}
	As in the proof of \Cref{lem:YQZNcommute}, the last term vanishes, and $\d f_\h(Y_j)=0$ for every $j$. Thus, the expression simplifies to
	\begin{equation*}
		[Y^\h_j,Y^\h_k]=\sum_l c^l_{jk}Y_l+\frac{1}{f_\h}(Y_k(f_{Y_j})-Y_j(f_{Y_k})) Z
		=\sum_l c^l_{jk}Y_l-\frac{2}{f_\h}\omega_\h(Y_j,Y_k)Z
	\end{equation*}
	Now we use \Cref{lem:twisted-structures}:
	\begin{align*}
		[Y_j^\q,Y^\q_k]&\sim_{\mc H}[Y_j^\h,Y_k^\h]
		+\frac{1}{f_\h}\omega_\h(Y_j^\h,Y_k^\h) Z
		=[Y_j^\h,Y_k^\h]
		+\frac{1}{f_\h}\omega_\h(Y_j,Y_k) Z\\
		&=\sum_l c^l_{jk}Y_l-\frac{1}{f_\h}\omega_\h(Y_j,Y_k) Z
	\end{align*}
	Since $-\big(\frac{1}{f_\h}Z\big)'=Z_{\bar N}$, we have $[Y^\q_j,Y^\q_k]=\sum_l c^l_{jk}Y'_l+\omega_\h'(Y_j',Y_k')Z_{\bar N}$. Using the relation $Y^\q_j=(Y^\h_j)'=Y_j'-(f_{Y_j}f_\h^{-1}Z)'=Y_j'+f_{Y_j}'Z_{\bar N}$ can also write this as
	\begin{equation*}
		[Y^\q_j,Y^\q_k]=\sum_l c^l_{jk}Y^\q_l 
		+\big(\omega_\h'(Y'_j,Y'_k)-\sum_l c^l_{jk} f_{Y_l}'\big)Z_{\bar N}
	\end{equation*}
	which was to be shown. 
	
	To check that the combination $\omega_\h'(Y'_j,Y'_k)-\sum_l c^l_{jk} f_{Y_l}'$ is indeed constant, it suffices to show that its twist is. Thus, we compute the differential of $\omega_\h(Y_j,Y_k)-\sum_l c^l_{jk}f_{Y_l}$. Applying Cartan's formula multiple times and using $\d \omega_\h=0$ we write the first term as
	\begin{equation*}
		L_{Y_k}\iota_{Y_j}\omega_\h-\iota_{Y_k}L_{Y_j}\omega_\h
		=\iota_{Y_j}L_{Y_k}\omega_\h+\iota_{[Y_k,Y_j]}\omega_\h
		-\iota_{Y_k}L_{Y_j}\omega_\h
	\end{equation*}
	Both the first and last term vanish, since the vector fields $Y_j$ all preserve $\omega_\h$. Thus, we find
	\begin{equation*}
		\d\Big(\omega_\h(Y_j,Y_k)-\sum_l c^l_{jk}f_{Y_l}\Big)=\omega_\h([Y_k,Y_j],\cdot)
		+\sum_l c^l_{jk} \omega_\h(Y_l,\cdot )=0
	\end{equation*}
	since $[Y_k,Y_j]=\sum_l c^l_{kj}Y_l$.
\end{myproof}

The constants $A_{jk}$ define a two-cocycle $\alpha\coloneqq \frac{1}{2}A_{jk} Y^*_j\wedge Y^*_k\in Z^2(\mf g)$, where $\{Y^*_j\}$ denotes the dual basis of $\mf g^*$. It is well-known that the corresponding cohomology class $[\alpha]$ vanishes if and only if the extension is trivial. In other words, we have:

\begin{cor}
	Under the assumptions of \Cref{thm:structure}, if the class $[\alpha]\in H^2(\mf g)$ is trivial, there exists a redefinition $\bar Y_j^\q\coloneqq Y_j^\q+\beta_j Z_{\bar N}$ such that $\bar{\mf{g}}\coloneqq \Span\{\bar Y_j^\q\mid j=1,\dots,k\}\subset \mf g_\q$ is a subalgebra isomorphic to $\mf g$, and $\mf g_\q=\bar{\mf g}\oplus \R Z_{\bar N}$. \proofclear
\end{cor}

\begin{cor}
	If $\mf g\subset \aut(N,f_Z)$ is semi-simple, then it gives rise to an isomorphic subalgebra $\bar{\mf g}\subset \aut(\bar N,Z_{\bar N})$. 
\end{cor}
\begin{myproof}
	This follows from Whitehead's lemma, which asserts that $H^2(\mf g)=0$ for every semi-simple Lie algebra $\mf g$.
\end{myproof}

\clearpage

\subsection{Application to the \texorpdfstring{$c$}{c}-map}

Now we apply the above results to the supergravity $c$-map. We consider a PSK manifold $\bar M$, which, by definition, is the K\"ahler quotient $M\sslash S^1$ of a CASK manifold $(M,g_M,\omega_M,\nabla,\xi)$ by the circle action generated by $-J\xi$, where $\xi$ is the Euler field on $M$, and $J$ its complex structure.

Recall from \Cref{sec:prelims} that applying the supergravity $c$-map to $\bar M$ yields the same result as applying the rigid $c$-map to $M$, and composing it with the HK/QK correspondence. This fact allows us to study the supergravity $c$-map in the presence of symmetries through the HK/QK correspondence, which we discussed in the previous section. We will now verify that (infinitesimal) automorphisms of the PSK manifold $\bar M$ naturally lead, through the rigid $c$-map, to the set-up considered in \Cref{sec:HKQKaut}.

An automorphism of the PSK manifold $\bar M$ is, by definition, induced by an automorphism $\varphi$ of the CASK manifold $M$. Any diffeomorphism of $M$ admits a canonical lift to a diffeomorphism $\Phi$ of its cotangent bundle, using the pullback on one-forms: If $\alpha_p$ is a covector at $p\in M$, we define $\Phi((p,\alpha_p))=(\varphi(p),(\varphi^{-1})^*\alpha_p)$ (or in other words, $\Phi=(\varphi^{-1})^*$).

\begin{lem}\label{lem:liftedsymm}
	Let $M$ be a CASK manifold and $N=T^*M$ its cotangent bundle, endowed with the rigid $c$-map hyper-K\"ahler structure. If $\varphi$ is an automorphism of the CASK structure on $M$ and $\Phi$ its canonical lift to $N$, then $\Phi$ preserves the full hyper-K\"ahler structure of $N$, i.e.\ $\Phi^*g=g$ and $\Phi^*\omega_i=\omega_i$, $i=1,2,3$.
\end{lem}
\begin{myproof}
	Since $\varphi$ preserves the special K\"ahler connection, $\Phi$ preserves the splitting \eqref{eq:splitting}. Thus, we may use the expressions in \Cref{eq:cotangentHKstr}. But then, $g$, $I_1$, and $I_2$ are explicitly defined in terms of tensor fields preserved by $\varphi$, hence they are preserved by $\Phi$.
\end{myproof}

Now we consider the infinitesimal case. Assume that the initial CASK manifold $M$ (or equivalently the PSK manifold $\bar M$) admits a one-parameter group of automorphisms $\varphi_t$, generated by some vector field $X$. The one-parameter group $\Phi_t$ of canonical lifts generates a vector field $Y\in \mf X(N)$, which can be expressed as follows:

\begin{lem}\label{lem:canonliftformula}
	Let $(\nabla X)^*:T^*M\to T^*M$ be the adjoint of the endomorphism $\nabla X$ of $TM$, and $\eta\in \Gamma(T^{\mc V}N)$ the (fiberwise) Euler field on $T^*M=N$. Under the isomorphism $T^{\mc V}N\cong \pi^*(T^*M)$, we may view $(\nabla X)^*$ as an endomorphism of the vertical tangent bundle, and hence apply it to $\eta$. Then $Y=\tilde X-(\nabla X)^*(\eta)$, where $\tilde X$ is the $\nabla$-horizontal lift of $X$.
\end{lem}
\begin{myproof}
	Since the proposed formula for $Y$ certainly specifies a lift of $X$, we need only check that it is the canonical one. It is a basic fact from symplectic geometry that a diffeomorphism of $T^*M$ is the canonical lift of a diffeomorphism of $M$ if and only if it preserves the tautological one-form $\lambda$. Therefore, it suffices to check $L_{\tilde X-(\nabla X)^*(\eta)}\lambda=0$.
	
	Since $\nabla$ is flat, we may choose local $\nabla$-affine coordinates $(q^i)$ on $M$. With respect to local canonical coordinates $(q^i,p_j)$ on $T^*M$, we have $\eta=\sum p_i\pd{}{p_i}$ and consequently
	\begin{equation*}
		\tilde X-(\nabla X)^*(\eta)=\sum \bigg(X^i \pd{}{q^i}-\pd{X^i}{q^j}p_i \pd{}{p_j}\bigg)
	\end{equation*}
	With respect to any set of canonical coordinates, $\lambda=\sum p_i\d q^i$, so we compute
	\begin{equation*}
		L_{\tilde X-(\nabla X)^*(\eta)}\lambda=\sum \bigg( p_i \d X^i - \pd{X^i}{q^j}p_i \d q^j\bigg)=0
	\end{equation*}
	proving our claim.
\end{myproof}

\Cref{lem:liftedsymm} shows that the canonical lift $Y$ of an infinitesimal CASK automorphism $X$ preserves the full hyper-K\"ahler structure of $N=T^*M$. Hence, to prove that our discussion of automorphisms under the HK/QK correspondence applies, we need only verify that the lifted action is $\omega_1$-Hamiltonian and preserves $f_Z$, the Hamiltonian function of $Z=-\widetilde{J\xi}$ with respect to $\omega_1$.

\begin{lem}\label{lem:liftedaut}
	Let $\varphi_t$, $\Phi_t$, and $Y\in \mf X(N)$ be as above. Then
	\begin{numberedlist}
		\item The action generated by $Y$ is $\omega_1$-Hamiltonian as well as $\omega_3$-Hamiltonian.
		\item $Y$ preserves $f_Z$.
	\end{numberedlist}
\end{lem}
\begin{myproof}\leavevmode
	\begin{numberedlist}
		\item Since $\omega_3$ is the canonical symplectic structure on $T^*M$, it is given by $\omega_3=-\d\lambda$. Since $Y$ is the canonical lift of $X$, it preserves $\lambda$, and is $\omega_3$-Hamiltonian with Hamiltonian function $-\lambda(Y)$.
		
		Since the special K\"ahler connection on $M$ is symplectic, $\omega_M$ has constant coefficients $(\omega_M)_{ij}$ with respect to any choice of $\nabla$-affine coordinates $(q^i)$. With respect to the corresponding canonical coordinates $(q^i,p_j)$ on $N$, we therefore have
		\begin{equation*}
			\omega_1=\sum \Big((\omega_M)_{ij}\d q^i\wedge\d q^j+(\omega_M)^{ij}\d p_i\wedge \d p_j\Big)
		\end{equation*}
		where $(\omega_M)^{ij}$ are the coefficients of the inverse matrix. Using \Cref{lem:canonliftformula}, we find
		\begin{equation*}
			\iota_Y\omega_1=\sum \bigg(2(\omega_M)_{ij}X^i \d q^j 
			-2(\omega_M)^{ij}\frac{\partial X^k}{\partial q^i} p_k \d p_j\bigg)
		\end{equation*}
		By \Cref{lem:CASKHam}, the first term is given by $-\pi^*\d f_X$, where $f_X$ is a Hamiltonian function for $X$ with respect to $\omega_M$. Since $X$ preserves the special K\"ahler connection, its coefficients with respect to the $\nabla$-affine coordinates are affine functions, i.e.\ $\pd{X^i}{q^j}$ is constant for every $i,j$. Therefore, the second term is given by $-\d\big(\sum S^{jk}p_kp_j\big)$, where we defined a new tensor field $S\in \Sym^2((T^{\mc V}N)^*)$ by setting $S^{jk}\coloneqq \sum(\omega_M)^{ij}\pd{X^k}{q^i}$. In invariant terms, $S=\omega_1((\nabla X)^*\cdot,\cdot )\big|_{T^{\mc V}N}$, where we regard $(\nabla X)^*$ as an endomorphism of $T^\mc VN$; its symmetry follows from the fact that $0=L_X\omega_M=(\nabla_X-\nabla X)\omega_M$ and that $\omega_M$ is parallel with respect to $\nabla$. Then $\sum S^{ij}p_ip_j=S(\eta,\eta)$, where $\eta=\sum p_i \pd{}{p_i}\in \Gamma(T^{\mc V}N)$ is the fiberwise Euler field. Thus, we have found
		\begin{equation*}
			\iota_Y\omega_1=-\d(\pi^*f_X+S(\eta,\eta)
		\end{equation*}
		In particular, $Y$ is $\omega_1$-Hamiltonian.
		\item Since $f_Z=-\frac{1}{2}\pi^*g_M(\xi,\xi)$ and $X$ preserves both $g_M$ and $\xi$, its canonical lift $Y$ automatically preserves $f_Z$.
	\end{numberedlist}
\end{myproof}

\begin{prop}\label{prop:lift}
	Given a PSK manifold $\bar M$ and underlying CASK manifold $M$, there is an injective linear map $\aut(\bar M)\to \aut(N,f_Z)$, where $N$ is the image of $M$ under the rigid $c$-map.
\end{prop}
\begin{myproof}
	By \Cref{lem:infPSKCASKaut}, there exists an injective linear map $\aut(\bar M)\to \aut(M)$. In case the special K\"ahler connection $\nabla$ on $M$ does not equal the Levi-Civita connection, this is even a canonical isomorphism. By \Cref{lem:liftedsymm,lem:liftedaut}, the canonical lift of any element of $\aut (M)$ to $N=T^*M$ lies in $\aut(N,f_Z)$. By functoriality of the pullback this defines a homomorphism $\aut(M)\to \aut(N,f_Z)$, which is obviously injective as well.
\end{myproof}

Combining this with our results from \Cref{sec:HKQKaut}, we have several interesting consequences:

\begin{thm}\label{thm:cmapsymm}
	Let $\bar M$ be any PSK manifold, and $\bar N$ its image under the supergravity $c$-map. Then there exists an injective, linear map from $\aut(\bar M)$ into $\aut(\bar N,Z_{\bar N})$, the algebra of Killing fields on $\bar N$ which commute with $Z_{\bar N}$. This is true for the quaternionic K\"ahler metrics on $\bar N$ induced by the (undeformed) supergravity $c$-map, as well as their one-loop deformations.
\end{thm}
\begin{myproof}
	This follows immediately by composing the maps constructed in \Cref{prop:lift} and \Cref{prop:autmaps} (note that the latter introduces a certain freedom of choice, cf.~\Cref{rem:notunique}).
\end{myproof}

\begin{thm}\label{thm:dimisom}
	Let $\bar M$ be a PSK manifold of real dimension $2n$ and set $m=\dim\aut (\bar M)$. Let $(\bar N,g_\text{FS}^c)$, $c\geq 0$, be its image under the (one-loop corrected) supergravity $c$-map. Then the isometry group of $(\bar N,g^c_\text{FS})$ has dimension at least $m+2n+3$.
\end{thm}
\begin{myproof}
	Thinking of $\bar N$ as a bundle over $\bar M$, we already know from \Cref{lem:fiberisometries} that, for every $c\geq 0$, we have an isometric action of $\Heis_{2n+3}$ (possibly with discrete stabilizer), which restricts trivially to $\bar M$. By contrast, all isometries generated by the Killing fields constructed by means of \Cref{thm:cmapsymm} certainly act non-trivially on the base PSK manifold, hence they are independent of such fiberwise isometries.
\end{myproof}

In case $c=0$ the group of fiberwise isometries is even of dimension $2n+4$ (again, see \Cref{lem:fiberisometries}), so we recover Proposition 23 of \cite{CDJL2017}:

\begin{cor}
	Let $\bar M^{2n}$ be a PSK manifold with $m=\dim \aut(\bar M)$. Then the isometry group of $(\bar N,g^0_\text{FS})$, its image under the (undeformed) supergravity $c$-map, has dimension at least $m+2n+4$. \proofclear
\end{cor}

The next corollary generalizes Corollary 24 of \cite{CDJL2017}, where only the case $c=0$ was considered:

\begin{cor}\label{cor:cohomogeneity}
	Let $\bar M^{2n}$ be a PSK manifold such that its automorphism group acts with cohomogeneity $k$. Then the isometry group of $(\bar N,g^c_\text{FS})$ acts with cohomogeneity at most $k+1$. If $c=0$, the isometry group acts with cohomogeneity at most $k$. \proofclear
\end{cor}

Recall that all known homogeneous quaternionic K\"ahler manifolds of negative scalar curvature are Alekseevsky spaces, that is, quaternionic K\"ahler manifolds which admit a simply transitive, completely solvable group of isometries \cite{Ale1975,dv1992,Cor1996}. It is known that, with exception of the quaternionic hyperbolic spaces, they can all be obtained from the $c$-map \cite{dv1992}.

\begin{cor}
	For all Alekseevsky spaces, with exception of the quaternionic hyperbolic spaces, the one-loop deformation of the Ferrara-Sabharwal metric defines a one-parameter deformation by complete quaternionic K\"ahler manifolds $(\bar N,g^c_\text{FS})$, $c\geq 0$, which admit a group of isometries acting with cohomogeneity one. \proofclear
\end{cor}

Since there are many examples of complete, simply connected PSK manifolds \cite{CHM2012}, including homogeneous spaces \cite{dv1992,AC2000} and manifolds of cohomogeneity one \cite{CDJL2017}, this result yields many examples of quaternionic K\"ahler manifolds with non-trivial fundamental group.

In the next section, we study the one-loop deformations of the non-compact Wolf spaces $SU(2,n)/S(U(2)\times U(n))$ in more detail. For this series of examples, we show that the above results are sharp, i.e.~that the one-loop deformed Ferrara--Sabharwal metric is of cohomogeneity precisely one.

\clearpage

\section{Examples}
\label{sec:ex}

Certainly the most basic example of a PSK manifold is the space consisting of a single point. Applying the $c$-map to it, we obtain the simplest and (certainly in the physics literature) most prominent example of the $c$-map. The resulting four-dimensional quaternionic K\"ahler manifold%
\footnote{As is usual, we call a four-manifold quaternionic K\"ahler if it is (anti-)self-dual and Einstein.}
is known as the universal hypermultiplet. The universal hypermultiplet, which we will denote by $\bar N_0$, is diffeomorphic to $\R^4$, but the (undeformed) Ferrara--Sabharwal metric casts it as a non-compact Wolf space, namely the complex hyperbolic plane $\C H^2=SU(1,2)/U(2)$. The one-loop deformed universal hypermultiplet is not symmetric, and in fact not even locally homogeneous. This follows from the results of \cite{CS2018}, which we will review below. Thus, even the trivial PSK manifold yields an interesting quaternionic K\"ahler manifold.

The (undeformed) universal hypermultiplet $\bar N_0=SU(1,2)/S(U(1)\times U(2))$ is the first in an infinite series of non-compact Wolf spaces. Indeed, for any $k\geq 0$ the $(4k+4)$-dimensional symmetric space
\begin{equation*}
	\bar N_k\coloneqq \frac{SU(k+1,2)}{S(U(k+1)\times U(2))}
\end{equation*}
is quaternionic K\"ahler. These spaces can be characterized as the only Riemannian manifolds which are simultaneously quaternionic K\"ahler of negative scalar curvature, and K\"ahler. 

Each $\bar N_k$ can be constructed by applying the $c$-map. Consider
\begin{equation*}
	M_k=\bigg\{(z_0,\dots,z_k)\in \C^{k+1}\ \bigg|\ \abs{z_0}^2>\sum_{j=1}^k \abs{z_j}^2\bigg\}
\end{equation*}
We can view this domain in $\C^{k+1}$ as an affine special K\"ahler manifold by equipping it with the trivial special K\"ahler connection $\nabla=\d$ and the standard flat, indefinite K\"ahler structure of signature $(k,1)$:
\begin{equation}\label{eq:exCASKstr}
\begin{aligned}
	g_{M_k}&=\sum_{j=1}^k\abs{\d z_j}^2-\abs{\d z_0}^2\\
	\omega_{M_k}&=\frac{i}{2}\bigg(\sum_{j=1}^k \d z_j\wedge \d \bar z_j-\d z_0 \wedge \d \bar z_0\bigg)		
\end{aligned}
\end{equation}
With the position vector field $\xi=\sum_{\mu=0}^k \Big(z_\mu \pd{}{z_\mu}+\bar z_\mu \pd{}{\bar z_\mu}\Big)$ as Euler field, $M_k$ is a CASK domain. The vector field $-J\xi$ is (up to sign) just the generator of the standard $U(1)$-action of the unit complex numbers on $\C^{k+1}$ which leaves the full CASK structure invariant, and by taking the K\"ahler quotient we obtain (the projective model of) complex hyperbolic space $\C H^k$, which is therefore a PSK domain. Applying the supergravity $c$-map to $\C H^k$, one obtains $\bar N_k$ (see, for instance, \cite{CDS2017}). In the following, our main interest is in the one-loop deformations of these symmetric metrics.

The coordinate $z_0$ on $M_k$ plays a distinguished role in the $c$-map construction. Accordingly, the case $k=0$ is qualitatively different and fundamentally simpler than its higher-dimensional generalizations. In this case, which we discuss first, we are able to compute the full isometry group of the one-loop deformed Ferrara--Sabharwal metric. Its action on $\bar N_0$ is of cohomogeneity one, and it follows from \Cref{cor:cohomogeneity} that the isometry group acts with cohomogeneity at most one on the higher-dimensional generalizations as well. In fact, our results in \cite{CST2020a} imply that it acts with cohomogeneity precisely one. We summarize these results below.

\subsection{The universal hypermultiplet}

Since the universal hypermultiplet is diffeomorphic to $\R_{>0}\times \R^3$, we may use global coordinates $(\rho,\tilde\phi,\zeta,\tilde\zeta)$, with respect to which the (one-loop deformed) Ferrara--Sabharwal metric takes the following form~\cite{CS2018}:
\begin{equation*}
	g^c_{\bar N_0}=\frac{1}{2\rho^2}\bigg[\frac{\rho+2c}{\rho+c}\d \rho^2
	+\frac{\rho+c}{\rho+2c}(\d \tilde \phi+\zeta \d \tilde \zeta -\tilde \zeta \d \zeta)^2
	+2(\rho+2c)\big(\d \tilde\zeta^2+\d \zeta^2\big)\bigg]
\end{equation*}
For $c=0$, this is the complex hyperbolic metric on $\C H^2$, normalized such that the reduced scalar curvature $\nu = \frac{1}{4n(n+2)}\operatorname{scal}$ (where $n$ is the quaternionic dimension) equals $-1$. On the other hand, as $c\to \infty$, $g^c_{\bar N_0}$ tends to a metric of constant curvature. While both of these limiting cases are well understood, less is known for finite $c>0$. In this section, we compute the full isometry group of this metric. First, let us note that, by \Cref{lem:fiberisometries}, we have a free and isometric action of the three-dimensional Heisenberg group, which is transitive on the level sets of the coordinate $\rho$.

\begin{lem}
	The norm of the curvature endomorphism $\mc R:\bigwedge^2T^*M\to \bigwedge^2 T^*M$ of $(\bar N_0,g^c_{\bar N_0})$ is given by
	\begin{equation*}
	\norm{\mc R}^2=6\bigg(1+\frac{\rho^6}{(\rho+2c)^6}\bigg)
	\end{equation*}
	In particular, for any $c>0$, $\norm{\mc R}^2$ is an injective function of $\rho$.
\end{lem}
\begin{myproof}
	An explicit eigenbasis for the (self-adjoint) curvature operator was constructed in~\cite{CS2018}. With our normalization, the corresponding eigenvalues are
	\begin{equation*}
		\lambda_1=-\bigg(1+2\frac{\rho^3}{(\rho+2c)^3}\bigg)\ ,\hspace{1.25cm}
		\lambda_2=-1\ ,\hspace{1.25cm}
		\lambda_3=-\bigg(1-\frac{\rho^3}{(\rho+2c)^3}\bigg)
	\end{equation*}
	with multiplicities $\mu_1=1$, $\mu_2=3$ and $\mu_3=2$, respectively. This means that
	\begin{equation*}
		\norm{\mc R}^2=\sum_i \mu_i\lambda_i^2=6\bigg(1+\frac{\rho^6}{(\rho+2c)^6}\bigg)
	\end{equation*}
	as claimed.
\end{myproof}

Since all scalar curvature invariants are invariant under isometries, we deduce:

\begin{cor}
	There are no isometries of the deformed Ferrara--Sabharwal metric which change the value of the coordinate $\rho$. In particular, $(\bar N_0,g_{\bar N_0}^c)$ is not a homogeneous space. \proofclear
\end{cor}

\begin{lem}
	Any isometry $\varphi$ of the one-loop deformed Ferrara--Sabharwal metric $g_{\bar N_0}^{c>0}$ is of the form $\varphi(\rho,\tilde\phi,\tilde\zeta,\zeta)=(\rho,\psi(\tilde\phi,\tilde\zeta,\zeta))$, where $\psi:\R^3\to \R^3$ preserves the symmetric bilinear form $\alpha_\lambda\coloneqq \lambda(\d\tilde\phi+\zeta\d \tilde\zeta-\tilde\zeta\d \zeta)^2+\d \tilde \zeta^2+\d \zeta^2$ for every $\lambda\in (0,c/8)$.
\end{lem}
\begin{myproof}
	We rewrite the metric as 
	\begin{equation*}
		g^c_{\bar N_0}=\frac{1}{2\rho^2}\bigg(\frac{\rho+2c}{\rho+c}\d\rho^2+2(\rho+2c)\alpha_{\lambda(\rho)}\bigg)
	\end{equation*}
	where $\lambda(\rho)=\frac{1}{2}\frac{\rho+c}{(\rho+2c)^2}\in (0,c/8)$ for all $\rho>0$. This makes it clear that any diffeomorphism of the given form is an isometry. 
	
	Conversely, any isometry $\varphi$ must preserve the level sets of $\rho$. Since it is an isometry, it furthermore preserves the unit normal bundle of each such hypersurface, or equivalently leaves the normal vector field $\pd{}{\rho}$ invariant up to sign. The first observation implies that $\varphi(\rho,\tilde\phi,\tilde\zeta,\zeta)=(\rho,\psi(\rho,\tilde\phi,\tilde\zeta,\zeta))$, where $\psi$ preserves $\alpha_\lambda$ for every $\lambda\in(0,c/8)$. The second fact implies that $\psi$ is independent of $\rho$, proving that $\varphi$ is of the claimed form.
\end{myproof}

\begin{lem}
	The Killing algebra of $(\bar N_0,g^{c>0}_{\bar N_0})$ is spanned by the Killing vectors 
	\begin{equation*}
		X_1\coloneqq \pd{}{\tilde\phi}\ , \hspace{0.75cm}
		X_2=\pd{}{\zeta}-\tilde\zeta \pd{}{\tilde \phi}\, \hspace{0.75cm}
		X_3=\pd{}{\tilde \zeta}+\zeta \pd{}{\tilde\phi}\ ,\hspace{0.75cm}
		X_4=\tilde\zeta\pd{}{\zeta}-\zeta\pd{}{\tilde\zeta}\\
	\end{equation*}
\end{lem}
\begin{myproof}
	We know that any Killing vector must be tangent to the level sets $\mc H_\rho$ of $\rho$. For any such hypersurface, $\{X_1,X_2,X_3\}$ provides a global trivialization of the tangent bundle. Therefore, any Killing vector is of the form $X=f_1 X_1+f_2 X_2 + f_3 X_3$ for $\rho$-independent functions $f_i\in C^\infty(\mc H_\rho)$, which we assume are not all constant. By the previous Lemma, $L_X\alpha_\lambda=0$. Combining this with the fact that each $X_i$, $i=1,2,3$, preserves $\alpha_\lambda$, we obtain
	\begin{equation*}
		2\lambda\big(\d f_1 -2\tilde\zeta \d f_2 +2 \zeta \d f_3\big)
		\big(\d \tilde\phi + \zeta\d\tilde\zeta - \tilde\zeta\d\zeta\big)
		+2\big(\d \tilde\zeta \vee \d f_3 +\d \zeta\vee \d f_2\big)=0
	\end{equation*}
	Since this equation must hold for every $\lambda\in(0,c/8)$, the two terms actually vanish independently:
	\begin{gather*}
		\big(\d f_1 -2\tilde\zeta \d f_2 +2 \zeta \d f_3\big)
		\big(\d \tilde\phi + \zeta\d\tilde\zeta - \tilde\zeta\d\zeta\big)=0\\
		\big(\d \tilde\zeta \vee \d f_3 +\d \zeta\vee \d f_2\big)=0
	\end{gather*} 
	The latter implies that $f_2$ depends only on $\tilde\zeta$ while $f_3$ only depends on $\zeta$, with the constraint that $\pd{f_2}{\tilde\zeta}=-\pd{f_3}{\zeta}$. This forces these derivatives to be constant, and we conclude that $f_2$ and $f_3$ are affine functions. By adding constant multiples of $X_2$ and $X_3$, we may in fact arrange that these functions are linear, i.e.~$f_2=k\tilde\zeta$ and $f_3=-k\zeta$ for some $k\in \R$. Plugging this into the first equation, we find
	\begin{equation*}
		\big(\d f_1 -k\d (\tilde\zeta^2+\zeta^2)\big)\big(\d \tilde\phi+\zeta \d \tilde\zeta-\tilde\zeta\d\zeta\big)=0
	\end{equation*}
	and deduce that $f_1=k (\tilde\zeta^2+\zeta^2)$, up to a constant which we may remove by adding a multiple of $X_1$. In conclusion, we have
	\begin{equation*}
		X=k\big((\tilde \zeta^2+\zeta^2)X_1 + \tilde\zeta X_2-\zeta X_3\big)=kX_4
	\end{equation*}
	This proves the claim.
\end{myproof}

\begin{thm}
	The isometry group of the one-loop deformed Ferrara--Sabharwal metric on the universal hypermultiplet is $\Heis_3\rtimes O(2)$.
\end{thm}
\begin{myproof}
	The isometry group acts on the Killing algebra by isomorphisms, and in particular preserves its center, which is spanned by $\pd{}{\tilde\phi}$. Thus, any isometry $\varphi$ of $g^{c>0}_{\bar N_0}$ must be of the form $\varphi(\rho,\tilde\phi,\tilde\zeta,\zeta)=(\rho,\chi(\tilde\phi,\tilde\zeta,\zeta),\psi(\tilde\zeta,\zeta))$, where $\chi:\R^3\to \R$ and $\psi:\R^2\to \R^2$ are smooth maps. By the same argument as in the proof of the preceding Lemma, $\varphi$ must also preserve the symmetric bilinear forms $(\d \tilde\phi+\zeta\d\tilde\zeta-\tilde\zeta \d\zeta)^2$ and $\d \tilde\zeta^2+\d \zeta^2$. Recognizing the latter as the Euclidean metric on $\R^2$, we see that $\psi$ must be a Euclidean motion. Thinking of $\R^2$ as $\C$, with complex coordinate $\xi=\tilde\zeta+i\zeta$, we have $\psi(\xi)=e^{i\theta}(\xi+u)$ or $\psi(\xi)=e^{-i\theta}(\bar \xi+\bar u)$, where $u\in \C$. 
	
	Imposing that $\varphi$ also preserves $(\d \tilde\phi+\Im(\bar \xi\d \xi))^2$, we find
	\begin{equation*}
		\d\big(\tilde\phi \pm  \chi(\tilde\phi,\xi)\big) = \Im (\xi\d\bar \xi)\pm \Im(\psi(\xi)\d \overline{\psi(\xi)})
	\end{equation*}
	We look for solutions of this equation, given our general form for $\psi(\xi)$. Since the left-hand side is exact, so must the right-hand side. For $\psi(\xi)=e^{i\theta}(\xi+u)$, this only happens if the negative sign is chosen, while for $\psi(\xi)=e^{-i\theta}(\bar\xi+\bar u)$ the positive sign must be chosen. This leads to the differential equation
	\begin{equation*}
		\d \big(\tilde\phi + \Im(u \bar \xi) \pm \chi(\tilde\phi,\xi)\big)=0
	\end{equation*}
	whose general solutions are easily read off:
	\begin{equation*}
		\chi(\tilde\phi,\xi)=\mp \big(\tilde\phi+ \Im(u \bar \xi)+k\big)
	\end{equation*}
	where $k\in \R$ is an arbitrary constant. The translations in $\tilde\phi$ and $\xi$ determine a normal subgroup isomorphic to $\Heis_3$, while rotations and reflections in the $(\tilde\zeta,\zeta)$-plane give rise to the $O(2)$-subgroup.
\end{myproof}

\subsection{Higher-dimensional hypermultiplet manifolds}

To obtain the higher-dimensional symmetric spaces $\bar N_k=\frac{SU(k+1,2)}{S(U(k+1)\times U(2))}$ and their one-loop deformations, we start with complex hyperbolic space $\CH^k$. Introducing complex coordinates $\{X_j\}$ on the unit ball $B^{2n}\subset \C^n$, the complex hyperbolic metric can be expressed as:
\begin{equation*}
	g_{\CH^k}=\frac{1}{1-\abs{X}^2}\bigg(\sum_j \abs{\d X_j}^2+\frac{1}{1-\abs{X}^2}\bigg|\sum_j\bar X_j \d X_j\bigg|^2\bigg)
\end{equation*}
where $\abs{X}=\sum_j \abs{X_j}^2<1$. On the corresponding CASK manifold $M_k$ we have the standard linear action of $SU(k,1)$. It descends to $\CH^k$, casting it as the Hermitian symmetric space $SU(k,1)/U(k)$. In fact, $SU(k,1)$ acts by automorphisms of the PSK structure because it preserves the full CASK structure on $M_k$ (here we use that, in this case, the special K\"ahler and Levi-Civita connection coincide), and so we see that $\CH^k$ is homogeneous as a PSK manifold.

The one-loop deformed Ferrara--Sabharwal metric on $\bar N_k\cong \CH^k\times \R_{>0}\times \R\times \R^{k+1}\times \R^{k+1}$ is given in global coordinates $(X_j,\rho, \tilde\phi,\zeta^I,\tilde\zeta_I)$ by the following bulky expression:
\begin{align*}
	g_{\bar{N}_k}^c&=\frac{\rho+c}{\rho}g_{\CH^k}+\frac{1}{4\rho^2}\frac{\rho+2c}{\rho+c}\d \rho^2\\
	&\quad +\frac{1}{4\rho^2}\frac{\rho+c}{\rho+2c}
	\bigg(\d \tilde\phi +\sum_{I=0}^k (\zeta^I \d \tilde\zeta_I-\tilde\zeta_I \d \zeta^I)+\frac{2c}{1-\abs{X}^2}\Im \Big[\sum_{j=1}^k \bar X_j \d X_j\Big]\bigg)^2\\
	&\quad +\frac{1}{2\rho}\bigg(\sum_{j=1}^k \big((\d \tilde \zeta_j)^2+(\d \zeta^j)^2\big) - (\d \tilde \zeta_0)^2-(\d \zeta^0)^2\bigg)\\
	&\quad +\frac{\rho+c}{\rho^2}\frac{1}{1-\abs{X}^2}
	\bigg|\d \tilde\zeta_0+i\d \zeta^0 +\sum_{j=1}^k X^j\big(\d \tilde \zeta_j-i\d \zeta^j\big)\bigg|^2
\end{align*}
\Cref{thm:cmapsymm} implies that the full group $SU(k,1)$ of PSK symmetries of $\CH^k$ lifts to a subgroup of the isometry group of this metric; its orbits are of course $2k$-dimensional. Adding to this the fiberwise isometries (cf.~\Cref{lem:fiberisometries}), we have found a group of isometries whose orbits are of dimension $4k+3$, i.e.~hypersurfaces, in agreement with the general result \Cref{cor:cohomogeneity}. 

As in the case of the universal hypermultiplet, we have constructed a group acting by isometries, we have found that the isometry group acts with cohomogeneity at most one. A natural question is whether this result can be improved, or in other words, whether it is possible to find additional isometries which act non-trivially on the coordinate $\rho$. The following Theorem, the proof of which is given in~\cite{CST2020a}, provides a negative answer in all dimensions:

\begin{thm}
	The  norm of the curvature operator $\mc R: \bigwedge^2T^*M\to \bigwedge^2 T^*M$ of $(\bar N_k,g^c_{\bar N_k})$ is given by
	\begin{equation*}
		\norm{\mc R}^2=n(5n+1)+3\bigg((n-1)\frac{\rho}{\rho+2c}+\frac{\rho^3}{(\rho+2c)^3}\bigg)^2
		+3\bigg((n-1)\frac{\rho^2}{(\rho+2c)^2}+\frac{\rho^6}{(\rho+2c)^6}\bigg)
	\end{equation*}
	\proofclear
\end{thm}

\begin{cor}
	For any $c>0$, the norm of the curvature operator of $(\bar N_k,g_{\bar N_k}^c)$ is an injective function of $\rho$. As a consequence, any isometry of $g^{c>0}_{\bar N_k}$ must preserve $\rho$.
\end{cor}
\begin{myproof}
	Up to an additive constant, the function $\abs{\mc R}^2:\R_{>0}\to \R_{>0}$ is the composition of two functions $\varphi,\psi:\R_{>0}\to \R_{>0}$, given by 
	\begin{gather*}
		\varphi(x)=\frac{x}{x+2c}\\
		\psi(x)=3\big((x^3+(n-1)x)^2+x^6+(n-1)x^2\big)
	\end{gather*}
	The first function is easily checked to be injective. For $\psi$, assume that $\psi(x_1)-\psi(x_2)=0$ for some $x_1,x_2\in \R_{>0}$. Then
	\begin{equation*}
		(x_1^3+(n-1)x_1)^2+x_1^6+(n-1)x_1^2-(x_2^3+(n-1)x_2)^2-x_2^6-(n-1)x_2^2=0
	\end{equation*}
	Clearly, this expression vanishes if $x_1=x_2$, and since $\psi$ is even, it also vanishes if $x_1=-x_2$. Factorizing, we find:
	\begin{equation*}
		(x_1-x_2)(x_1+x_2)\Big(2(x_1^4+x_2^4)+2x_1^2x_2^2+2(x_1^2+x_2^2)(n-1)+n(n-1)\Big)=0
	\end{equation*}
	The second and third factor are each manifestly positive for $x_1,x_2\in \R_{>0}$, so the only possibility is $x_1=x_2$.
\end{myproof}

As an immediate consequence, we have:

\begin{thm}
	The one-loop deformation $g^{c>0}_{\bar N_k}$ of the Ferrara--Sabharwal metric on $\bar N_k$ is of cohomogeneity one. In particular, $(\bar N_k,g^{c>0}_{\bar N_k})$ is not a homogeneous space.
\end{thm}

\end{document}